\documentclass[11pt]{article}
\usepackage{amsfonts}

\usepackage{amsmath}
\usepackage{latexsym,amsfonts,amssymb}

\newcommand{\epsln}{\varepsilon}

\newcommand{\bfa}{{\mathbf a}}
\newcommand{\bfb}{{\mathbf b}}
\newcommand{\bfe}{{\mathbf e}}
\newcommand{\bfu}{{\mathbf u}}
\newcommand{\bfv}{{\mathbf v}}
\newcommand{\bfw}{{\mathbf w}}
\newcommand{\bfg}{{\mathbf g}}
\newcommand{\bfh}{{\mathbf h}}

\newcommand{\bfn}{{\mathbf n}}
\newcommand{\bfI}{{\mathbf I}}

\newcommand{\bfT}{{\mathbf T}}

\newcommand{\hs}{\hspace{0.3cm}}

\newcommand{\vsss}{\vspace{.05in}}

\font\tensy=cmsy10

\begin{document}

\title{Asymptotic Behavior of Solutions of a Free Boundary Problem
  Modelling the Growth of Tumors with Stokes Equations}
\author{Junde Wu$^{\dag}$\ \ \ and\ \ \ Shangbin Cui$^{\ddag}$\\[0.2cm]
  {\small $^{\dag}$ Department of Mathematics, Sun Yat-Sen University,
  Guangzhou, Guangdong 510275,}\\ [-0.1cm]
  {\small People's Republic of China. E-mail:\,wjdmath@yahoo.com.cn}
  \\[0.1cm]
  {\small $^{\ddag}$ Institute of Mathematics, Sun Yat-Sen University,
  Guangzhou, Guangdong 510275,}\\ [-0.1cm]
  {\small People's Republic of China. E-mail:\,cuisb3@yahoo.com.cn}}
\date{}
 \maketitle

\begin{abstract}
   We study a free boundary problem modelling the growth of non-necrotic
   tumors with fluid-like tissues. The fluid velocity satisfies Stokes
   equations with a source determined by the proliferation rate of tumor
   cells which depends on the concentration of nutrients, subject to a
   boundary condition with stress tensor effected by surface tension. It
   is easy to prove that this problem has a unique radially symmetric
   stationary solution. By using a functional approach, we prove that there
   exists a threshold value $\gamma_*>0$ for the surface tension coefficient
   $\gamma$, such that in the case $\gamma>\gamma_*$ this radially symmetric
   stationary solution is asymptotically stable under small non-radial
   perturbations, whereas in the opposite case it is unstable.
\medskip

\medskip


\end{abstract}

\section{Introduction}
\setcounter{equation}{0}

  In this paper we study the following free boundary problem modelling
  the growth of non-necrotic tumors with fluid-like tissues:
\begin{equation}
   \Delta\sigma=f(\sigma) \quad\mbox{in}\;\;
    \Omega(t),\;\;t>0,
\end{equation}
%(1.1)
\begin{equation}
  \nabla\cdot\bfv=g(\sigma) \quad\mbox{in}\;\;
  \Omega(t),\;\;t>0,
\end{equation}
%(1.2)
\begin{equation}
   -\nu\Delta\bfv+\nabla p-{\nu\over3}\nabla
   (\nabla\cdot\bfv)=0 \quad\mbox{in} \;\;\Omega(t),\;\;t>0,
\end{equation}
%(1.3)
\begin{equation}
   \sigma=\bar\sigma \quad\mbox{on}\;\;\partial\Omega(t),\;\;t>0,
\end{equation}
%(1.4)
\begin{equation}
  \bfT({\bfv},p){\bfn}=-\gamma\kappa{\bfn}
  \quad\mbox{on}\;\;\partial\Omega(t),\;\;t>0,
\end{equation}
%(1.5)
\begin{equation}
   V_n=\bfv\cdot\bfn\quad\mbox{on}\;\;\partial\Omega(t),\;\;t>0,
\end{equation}
%(1.6)
\begin{equation}
   \int_{\Omega(t)}\bfv\;dx=0, \quad t>0,
\end{equation}
%(1.7)
\begin{equation}
   \int_{\Omega(t)}\bfv\times{x}d\;x=0,\quad t>0,
\end{equation}
%(1.8)
\begin{equation}
   \Omega(0)=\Omega_0,
\end{equation}
%(1.9)
  where $\sigma=\sigma(t,x)$, $\bfv=\bfv(t,x)$ ($=(v_1(t,x),v_2(t,x),v_3(t,x))$)
  and $p=p(t,x)$ are unknown functions representing the concentration of nutrient,
  the velocity of the fluid and the internal pressure, respectively, $f$ and $g$
  are given functions representing the nutrient consumption rate and tumor cell
  proliferation rate, respectively, which typically have the following forms
  respectively:
\begin{equation}
   f(\sigma)=\lambda\sigma, \quad
   g(\sigma)=\mu(\sigma-\sigma_c),
\end{equation}
%(1.10)
  where $\lambda$, $\mu$ and $\sigma_c$ are positive constants, $\sigma_c<\bar\sigma$,
  and $\Omega(t)$ is an a priori unknown bounded domain in $\Bbb R^3$ representing
  the region occupied by the tumor at time $t$. Besides, $\nu$, $\bar\sigma$ and
  $\gamma$ are positive constants, among which $\nu$ is the viscosity coefficient
  of the fluid, $\gamma$ is the surface tension coefficient of the tumor surface,
  and $\bar\sigma$ is the concentration of nutrient in tumor's host tissues, $\kappa$,
  $V_n$ and $\bfn$ denote the mean curvature, the normal velocity and the unit outward
  normal, respectively, of the tumor surface $\partial\Omega(t)$, and $\bfT({\bfv},p)$
  represents the stress tensor, i.e.,
\begin{equation}
  \bfT({\bfv},p)=\nu\big[\nabla\otimes\bfv+(\nabla\otimes\bfv)^T\big]-
  (p+{2\nu\over3}\nabla\cdot\bfv)\bfI,
\end{equation}
%(1.11)
  where $\bfI$ denotes the unit tensor. We note that the sign of the
  mean curvature $\kappa$ is defined such that it is nonnegative for
  convex hyper-surfaces. Without loss of generality, later on we assume
  that
$$
  \nu=1 \quad \mbox{and} \quad \bar{\sigma}=1.
$$
  Note that the general situation can be easily reduced into this special
  situation by using the rescaling $\sigma\to\sigma/\bar{\sigma}$, $p\to p/\nu$
  and $\gamma\to\gamma/\nu$.

  Tumor growth modelling and analysis has attracted considerable attention
  during the past more than ten years. Most tumor models assume that the
  tumor tissue has the structure of a porous medium for which Darcy's law
  applies (see, e.g., \cite{Byrne}, \cite{ByrnChap} and the references
  cited therein). For such tumor models, many interesting results of
  rigorous analysis have been obtained, for which we refer the interested
  reader to see \cite{Cui0}--\cite{CuiEsc2}, \cite{FriedHu1}--\cite{FriedReit3},
  \cite{WuCui}, \cite{ZhouCui} and the references cited therein. The tumor
  whose tissue does not have the structure of a porous medium but instead
  is more like a fluid was recently considered by Franks et al in the
  literatures \cite{Franks1}--\cite{Franks4}, where some new models
  were proposed to mimic the early stages of the growth of ductal
  carcinoma in the breast. A basic feature of a ductal carcinoma
  in the breast in early stages is that it is confined to the duct
  of a mammary gland, which consists of epithelial cells, a meshwork
  of proteins, and extracellular fluid. In modelling, this leads to
  the replacement of the Darcy's law used in porous medium structured
  tumor models by the Stokes equations. See \cite{Franks1}--\cite{Franks4}
  for details. The models of Franks et al \cite{Franks1}--\cite{Franks4}
  have been concisely reformulated by Friedman in \cite{Fried1} (see
  also \cite{Fried2}).

  The problem (1.1)--(1.9) above is a simplification of the tumor model proposed
  in \cite{Fried1}. The simplifications are made in two aspects. First, the model
  in \cite{Fried1} contains a system of nonlinear hyperbolic conservation laws
  with source terms describing the movements and interchanges of three different
  species of cells: proliferating cells, quiescent cells and necrotic cells. In
  this paper we only consider one species of proliferating cells, so that no
  hyperbolic conservation laws appear in (1.1)--(1.9). Second, in \cite{Fried1}
  the equation for $\sigma$ is of the following evolutionary type:
$$
   \partial_t\sigma=\Delta\sigma-f(\sigma) \quad\mbox{in}\;\;
   \Omega(t),\;\;t>0,
$$
  where $f$ is as given in (1.10), but in this paper the stationary form (1.1)
  is considered. All these simplifications are made for the purpose to make the
  analysis simpler. If either one of the above two aspects of simplifications
  are not made, then the model will be much more complicated to analyze, and
  new mathematical techniques have to be employed. We leave it for future work.

  In \cite{Fried1} Friedman established local wellposedness in H\"{o}lder
  spaces of his model. Meanwhile, he proved that in the special case that
  the tumor contains only one species of cells (i.e., the tumor contains
  only proliferating cells), there exists a unique radially symmetric
  stationary solution. Based on these results, a number of interesting
  questions are raised in \cite{Fried1} (see also \cite{Fried2}), one
  of which is as follows: Is this radially symmetric stationary solution
  asymptotically stable under non-radial perturbations? A heuristic result
  toward an answer to this question was obtained by Friedman and Hu in
  \cite{FriedHu2}, where they proved that this radially symmetric stationary
  solution is {\em linearly asymptotically stable} for small $\mu/\gamma$,
  i.e., there exists a threshold value $(\mu/\gamma)_*$ such that if we
  denote by $(\sigma_s,\bfv_s,p_s,\Omega_s)$ this stationary solution,
  then in the case $\mu/\gamma<(\mu/\gamma)_*$ the trivial solution of
  the linearization at $(\sigma_s,\bfv_s,p_s,\Omega_s)$ of the original
  problem is asymptotically stable. Moreover, they also proved that in
  the case $\mu/\gamma>(\mu/\gamma)_*$ the radially symmetric stationary
  solution is unstable. We also refer the interested reader to see
  \cite{FriedHu1} for the study of existence of non-radial stationary
  solutions.

  The purpose of this paper is to prove that, at least for the simplified
  model (1.1)--(1.9), the answer to the above question is yes for large
  $\gamma$ but no for small $\gamma$. To this end, we shall use a functional
  approach inherited from the references \cite{CuiEsc1}, \cite{EscSim},
  \cite{WuCui} and \cite{ZhouCui} to study the problem (1.1)--(1.9), namely,
  we shall first reduce the problem (1.1)--(1.9) into an evolution equation
  containing merely the function $\rho$ describing the free boundary
  $\partial\Omega(t)$,
  which can be considered as a differential equation in certain Banach space.
  We shall prove that this differential equation is of the parabolic type.
  Next we use the geometric theory for parabolic differential equations in
  Banach spaces (see \cite{Amann} and \cite{Lunar}) to study the stability
  of the stationary solution. Since our discussion does not depend on the
  specific linear forms of the equations (1.1) and (1.2), throughout this
  paper we shall not consider the specific forms of $f$ and $g$ given by
  (1.10), but instead assume that they are general smooth functions satisfying
  the following assumptions:

\medskip
  $(A1)$ $f\in C^\infty[0,\infty)$, $f'(\sigma)>0$ for $\sigma\ge0$
  and $f(0)=0$.

  $(A2)$ $g\in C^\infty[0,\infty)$, $g'(\sigma)>0$ for $\sigma\ge0$
  and $g(\sigma_c)=0$ for some $\sigma_c>0$,

  $(A3)$ $\sigma_c<1$.
\medskip

\noindent
  To give a precise statement of our main result, let us first introduce
  some notations.

  Given a bounded domain $\Omega\subseteq\Bbb R^3$ and two numbers
  $m\in \Bbb N$ and $\theta\in(0,1)$, we denote by $h^{m+\theta}
  (\overline\Omega)$ the so-called little H\"{o}lder space on
  $\Omega$ of index $m+\theta$, which is, by definition, the
  closure of $C^{\infty}(\overline\Omega)$ in the usual H\"{o}lder
  space $C^{m+\theta}(\overline\Omega)$. Similarly, given a smooth
  hypersurface $\Gamma$ in $\Bbb R^3$, we denote by $h^{m+\theta}
  (\Gamma)$ the closure of $C^\infty(\Gamma)$ in $C^{m+\theta}(\Gamma)$.

  It can be easily shown (see Theorem A in Appendix A) that under
  Assumptions $(A1)$--$(A3)$, the problem (1.1)--(1.8) has a unique
  radially symmetric stationary solution. Later on we use the same
  notation $(\sigma_s,\bfv_s,p_s,\Omega_s)$ as before to denote this
  radially symmetric stationary solution of the problem (1.1)--(1.8).
  Note that this means that there exists $R_s>0$ such that
  $\Omega_s=\{x\in\Bbb R^3:\,|x|<R_s\}$ and
\begin{equation}
  \sigma_s(x)=\sigma_s(r), \quad \bfv_s(x)=v_s(r)\frac{x}{r}, \quad
  p_s(x)=p_s(r) \quad \mbox{for}\;\; x\in\Omega_s,
\end{equation}
%(1.12)
  where $r=|x|$ and $v_s$ represents a scalar function. Clearly,
  a coordinate translation of a solution of (1.1)--(1.8) is still
  a solution of it. Thus, for any $x_0\in \Bbb R^3$, we denote by
  $(\sigma_{[x_0]},\bfv_{[x_0]},p_{[x_0]},\Omega_{[x_0]})$ the
  stationary solution obtained by the coordinate translation
  $x\to x+x_0$ of the stationary solution $(\sigma_s,\bfv_s,p_s,
  \Omega_s)$. Given $\rho\in C^1(\partial\Omega_s)$ with $\|\rho\|
  _{C^1(\partial\Omega_s)}$ sufficiently small, we denote by
  $\Omega_\rho$ the domain enclosed by the hypersurface $r=R_s+
  \rho(\omega)$, where $\omega\in\partial\Omega_s$. It is obvious that
  for $x_0\in \Bbb R^3$ such that $|x_0|$ is sufficiently small,
  there exists a smooth function $\rho_{[x_0]}$ on $\partial\Omega_s$
  such that $\Omega_{[x_0]}=\Omega_{\rho_{[x_0]}}$. Since we shall
  only be concerned with small perturbations of the stationary solution
  $(\sigma_s,\bfv_s,p_s,\Omega_s)$, it is natural to assume that the
  domains $\Omega(t)$ and $\Omega_0$ in (1.1)--(1.9) are small
  perturbations of $\Omega_s$. It follows that there exist functions
  $\rho(t)$ ($=\rho(\omega,t)$) and $\rho_0$ ($=\rho_0(\omega)$) on
  $\partial\Omega_s$ such that $\Omega(t)=\Omega_{\rho(t)}$ and
  $\Omega_0=\Omega_{\rho_0}$. Using these notations, the initial
  condition (1.9) can be rewritten as follows:
\begin{equation}
  \rho(\omega,0)=\rho_0(\omega) \quad \mbox{for}
  \;\; \omega\in\partial\Omega_s.
\end{equation}
%(1.13)
  The solution $(\sigma,{\bfv},p,\Omega)$ of the problem (1.1)--(1.9) will be
  correspondingly rewritten as $(\sigma,{\bfv},p,\rho)$, and the radially
  symmetric stationary solution $(\sigma_s,\bfv_s,p_s,\Omega_s)$ will be
  re-denoted as $(\sigma_s,\bfv_s,p_s,0)$.

  The  main result of this paper is the following theorem:
\medskip

  {\bf Theorem 1.1}\ \ {\em Assume that Assumptions $(A1)$--$(A3)$ hold. For
  given $m\in {\Bbb N}$, $m\geq 3$, and $0<\theta<1$, we have the following
  assertion: There exists a positive threshold value $\gamma_*$ such that for
  any $\gamma>\gamma_*$, the radially symmetric stationary solution $(\sigma_s,
  v_s,p_s,0)$ is asymptotically stable in the following sense: There exists
  constant $\varepsilon>0$ such that for any $\rho_0\in h^{m+\theta}
  (\partial\Omega_s)$ satisfying $\|\rho_0\|_{C^{m+\theta}(\partial\Omega_s)}
  <\varepsilon$, the problem $(1.1)$--$(1.9)$ has a unique solution $(\sigma,
  \bfv,p,\rho)$ for all $t\geq 0$, and there are positive constants $\omega$,
  $K$ independent of the initial data and a point $x_0\in\Bbb R^3$ uniquely
  determined by the initial data, such that the following holds for all $t\geq
  0$:
\begin{eqnarray}
  &&\|\sigma(\cdot,t)-\sigma_{[x_0]}\|_{C^{m+\theta}(\Omega(t))}+
  \|\bfv(\cdot,t)-\bfv_{[x_0]}\|_{C^{m-1+\theta}(\Omega(t))} \nonumber
  \\
  &+&\|p(\cdot,t)-p_{[x_0]}\|_{C^{m-2+\theta}(\Omega(t))}
  +\|\rho(\cdot,t)-\rho_{[x_0]}\|_{C^{m+\theta}(\partial\Omega_s)}
  \le K e^{-\omega t}.
\end{eqnarray}
%(1.14)
  For $\gamma<\gamma_*$ the stationary solution $(\sigma_s,{\bfv}_s,p_s,0)$
  is unstable.$\qquad$$\Box$}
\medskip

  It is interesting to compare this result with the corresponding result for
  the porous medium structured tumor model obtained by Cui and Escher in
  \cite{CuiEsc2}, where it is proved that, for the porous medium structured
  tumor model, there exists a threshold value for the surface tension
  coefficient $\gamma$, which we denote as $\tilde\gamma_*$, such that the
  unique radially symmetric stationary solution is asymptotically stable if
  $\gamma>\tilde\gamma_*$, but unstable if $\gamma<\tilde\gamma_*$. We shall
  show that $\gamma_*>\tilde\gamma_*$. This implies that radially symmetric
  stationary solution is more stable for a tumor whose tissue has a porous
  medium structure than a tumor whose tissue is more like a fluid. See Lemma
  3.5 for the proof of the assertion that $\gamma_*>\tilde\gamma_*$.

  The structure of the rest part is as follows. In Section 2 we first convert
  the problem into an equivalent initial-boundary value problem on a fixed
  domain by using the so-called Hanzawa transformation, and next further reduce
  it into a scalar equation containing the single function $\rho$, which can be
  regarded as a differential equation in the Banach space $h^{m-1+\theta}
  ({\Bbb S}^2)$. We shall also prove that this equation is of the parabolic
  type. In Section 3 we study the linearization of (1.1)--(1.8) at the radially
  symmetric stationary solution, and study the spectrum of the linearized
  operator. In the last section we give the proof of Theorem 1.1.

\section{Reduction of the problem}
\setcounter{equation}{0}

  In this section we reduce the problem (1.1)--(1.9) into a differential
  equation in a Banach space. For simplicity of the notation, later on
  we always assume that $R_s=1$. Note that this assumption is reasonable
  because the case $R_s\neq 1$ can be easily reduced into this case after
  a rescaling. It follows that
$$
  \Omega_s=\Bbb B^3=\{x\in\Bbb R^3:|x|<1\}\quad \mbox{and} \quad
  \partial\Omega_s=\Bbb S^2=\{x\in\Bbb R^3:|x|=1\}.
$$

  Let $m$ and $\theta$ be as in Theorem 1.1. For $u\in h^{m+\theta}
  (\overline{\Bbb B}^3)$, we denote by $\mbox{tr}_{\Bbb S^2}(u)$ the
  trace of $u$ on $\Bbb S^2$, i.e., $\mbox{tr}_{\Bbb S^2}(u)=u|_{\Bbb
  S^2}$. We know that $\mbox{tr}_{\Bbb S^2}(u)\in h^{m+\theta}(\Bbb S^2)$
  and the operator $\mbox{tr}_{\Bbb S^2}:u\to\mbox{tr}_{\Bbb S^2}(u)$
  from $h^{m+\theta}(\overline{\Bbb B}^3)$ to $h^{m+\theta}(\Bbb S^2)$
  is linear, bounded and  surjective. Let $\Pi$ be a bounded right
  inverse of it, i.e., $\Pi\in L(h^{m+\theta}(\Bbb S^2),h^{m+\theta}
  (\overline{\Bbb B}^3))$ and $\mbox{tr}_{\Bbb S^2}\big(\Pi(u))\big)=u$
  for any $u\in h^{m+\theta}(\Bbb S^2)$. Let $E\in L(C^{m+\theta}(\overline
  {\Bbb B}^3),BUC^{m+\theta}({\Bbb R}^3))$ be an extension operator, i.e.,
  $E$ has the property that $E(u)(x)=u(x)$ for any $u\in C^{m+\theta}
  (\overline{\Bbb B}^3)$ and $x\in\overline{\Bbb B}^3$. Here $BUC^{m+
  \theta}({\Bbb R}^3)$ denotes the space of all $C^{m}$ functions $u$
  on ${\Bbb R}^3$ such that $u$ itself and all its partial derivatives
  of order$\leq m$ are bounded and uniformly $\theta$-th order H\"{o}lder
  continuous in $\Bbb R^3$. We denote $\Pi_1=E\circ\Pi$. Then clearly
  $\Pi_1\in L(h^{m+\theta}(\Bbb S^2),h^{m+\theta}({\Bbb R}^3))$, where
  $h^{m+\theta}({\Bbb R}^3)$ represents the closure of $BUC^\infty
  ({\Bbb R}^3)$ in $BUC^{m+\theta}({\Bbb R}^3)$. Hence there exists a
  constant $C_0>0$ such that
\begin{equation}
  \|\Pi_1(\rho)\|_{BUC^{m+\theta}({\Bbb R}^3)}
  \le C_0\|\rho\|_{C^{m+\theta}(\Bbb S^2)}
  \quad\mbox{for}\;\;\rho\in h^{m+\theta}(\Bbb S^2).
\end{equation}
%(2.1)
  Take a constant $0<\delta<\min\{1/6,1/(3C_0)\}$ and fix it, where $C_0$ is
  the constant in (2.1). We choose a cut-off function $\chi\in C^\infty[0,
  \infty)$ such that
\begin{equation}
  0\le\chi\le1, \quad \chi(\tau)=\left\{
  \begin{array}{ll}
  1,\quad \text{for}\;|\tau|\le\delta,
  \\
  0,\quad \text{for}\;|\tau|\ge3\delta,
  \end{array}
  \right.
  \quad\text{and}\quad
  |\chi'(\tau)|\le{2\over3\delta}.
\end{equation}
%(2.2)
  We denote
$$
  O_\delta^{m+\theta}(\Bbb S^2)=\{\rho\in h^{m+\theta}(\Bbb S^2):
  \|\rho\|_{C^{m+\theta}({\Bbb S^2})}<\delta\}.
$$
  Given $\rho\in O_\delta^{m+\theta}(\Bbb S^2)$, we define the
  Hanzawa transformation $\Phi_\rho:{\Bbb R}^3\to\Bbb R^3$
  as follows:
\begin{equation}
  \Phi_\rho(x)=x+\chi(r-1)\Pi_1(\rho)(x)\frac{x}{r}
  \quad\mbox{for}\;\; x\in {\Bbb R}^3.
\end{equation}
%(2.3)
  Using (2.1) and (2.2) we can easily verify that $\Phi_\rho$ is a
  $h^{m+\theta}$ diffeomorphism from ${\Bbb R}^3$ onto itself, i.e.,
  $\Phi_\rho\in {\rm Diff}^{m+\theta}({\Bbb R}^3,{\Bbb R}^3)$, and
  each component of $\Phi_\rho$ and $\Phi_\rho^{-1}$ belongs to
  $h^{m+\theta}({\Bbb R}^3)$. Later on we write $\Phi_\rho\in
  \mbox{Diff}^{m+\theta}_h({\Bbb R}^3,{\Bbb R}^3)$ to indicate
  this fact. We define $\phi_\rho=\Phi_\rho\big|_{\Bbb S^2}$,
  and denote
$$
  \Omega_\rho=\Phi_\rho({\Bbb B}^3), \quad \Gamma_\rho=\partial
  \Omega_\rho=\Phi_\rho({\Bbb S}^2)=\mbox{Im}(\phi_\rho).
$$
  Clearly,
$$
  \phi_\rho(\omega)=[1+\rho(\omega)]\omega \quad
  \mbox{for}\;\; \omega\in\Bbb S^2.
$$
  This implies that $x\in\Gamma_\rho$ if and only if there exists
  $\omega\in\Bbb S^2$ such that $x=[1+\rho(\omega)]\omega$. Thus,
  in the polar coordinates $(r,\omega)$ of $\Bbb R^3$, where $r=|x|$
  and $\omega=x/|x|$, the hypersurface $\Gamma_\rho$ has the following
  equation: $r=1+\rho(\omega)$.

  Next, given $\rho\in C([0,T],O_\delta^{m+\theta}(\Bbb S^2))$, for
  each $t\in[0,T]$ we denote
$$
  \Gamma_\rho(t)=\Gamma_{\rho(t)} \quad \mbox{and} \quad
  \Omega_\rho(t)=\Omega_{\rho(t)}.
$$
  Since our purpose is to study asymptotical stability of the radially
  symmetric stationary solution, later on we always assume that the
  initial domain $\Omega_0$ is contained in a small neighborhood of
  $\Omega_s=\Bbb B^3$. It follows that there exists $\rho_0\in O_
  \delta^{m+\theta}(\Bbb S^2)$ such that $\Gamma_0\equiv\partial\Omega_0
  =\Gamma_{\rho_0}$.

  Let $\rho$ be as above, and let $\Phi_\rho^i$ be the $i$-th component
  of $\Phi_\rho$, $i=1,2,3$. We denote
$$
  [D\Phi_\rho]_{ij}:=\partial_i\Phi_\rho^j={\partial\Phi_\rho^j\over
  \partial x_i},\quad a_{ij}^\rho(x)=[D\Phi_\rho(x)]^{-1}_{ij}\quad
  (i,j=1,2,3),
$$
$$
  G_\rho(x)=\det(D\Phi_\rho(x)\big)\quad \mbox{for}\;\; x\in \Bbb R^3,
$$
$$
  H_\rho(\omega)=|\phi_\rho|^2\sqrt{1+|\nabla_\omega\phi_\rho|^2}
  \quad \mbox{for}\;\; \omega\in {\Bbb S}^2,
$$
  where $\nabla_\omega$ represents the orthogonal projection of the gradient
  $\nabla_x$ onto the tangent space $T_\omega({\Bbb S}^2)$\footnotemark[1]$^)$.
\footnotetext[1]{$^)$In the coordinate $\omega=\omega
  (\vartheta,\varphi)=(\sin\vartheta\cos\varphi,\sin\vartheta\sin\varphi,
  \cos\vartheta)$ ($0\leq\vartheta\leq\pi$, $0\leq\varphi\leq 2\pi$) of the
  sphere we have
$$
  \nabla_\omega f(\omega)=(\cos\vartheta\cos\varphi,\cos\vartheta
  \sin\varphi,-\sin\vartheta)\partial_\vartheta f(\omega(\vartheta,
  \varphi))+{1\over\sin\vartheta}(-\sin\varphi,\cos\varphi,0)
  \partial_\varphi f(\omega(\vartheta,\varphi)).
$$
Note that $\nabla_x f=\displaystyle\frac{\partial f} {\partial
r}\omega+{1\over r}\nabla_\omega f$.
}
  Here and hereafter, for a matrix $A$ we use the notation $A_{ij}$
  to denote the element of $A$ in the $(i,j)$-th position. Since
  $\Phi_\rho\in\mbox{Diff}^{m+\theta}_h({\Bbb R}^3,{\Bbb R}^3)$,
  we have $a_{ij}^\rho\in h^{m-1+\theta}({\Bbb R}^3)$, $i,j=1,2,3$,
  $G_\rho\in h^{m-1+\theta}({\Bbb R}^3)$, and $H_\rho\in h^{m-1+\theta}
  ({\Bbb S}^2)$. We now introduce four partial differential operators
  ${\mathcal A}(\rho)$, $\vec{\mathcal B}(\rho)$, $\vec{\mathcal B}(\rho)
  \cdot$ and $\vec{\mathcal B}(\rho)\otimes$ on ${\Bbb R}^3$ as follows:
$$
  {\mathcal A}(\rho)u(x)=a^\rho_{ij}(x){\partial_j}
  \big(a^\rho_{ik}(x){\partial_k u(x)}\big) \quad
  \mbox{for scalar function}\;\; u,
\eqno{(2.4)}
$$
$$
  \vec{\mathcal B}(\rho)u(x)=\big(a^\rho_{1j}(x)
  {\partial_j u(x)}, a^\rho_{2j}(x){\partial_j u(x)},
  a^\rho_{3j}(x){\partial_j u(x)}\big) \quad
  \mbox{for scalar function}\;\; u,
\eqno{(2.5)}
$$
$$
  \vec{\mathcal B}(\rho)\cdot{\bfv}(x)=a^\rho_{ij}(x)
  {\partial_j v_i(x)} \quad
  \mbox{for vector function}\;\; {\bfv}=(v_1,v_2,v_3).
\eqno{(2.6)}
$$
$$
  \vec{\mathcal B}(\rho)\otimes{\bfv}(x)=(a^\rho_{ik}(x)
  {\partial_k v_j(x)}) \quad
  \mbox{for vector function}\;\; {\bfv}=(v_1,v_2,v_3).
\eqno{(2.7)}
$$
  Here and hereafter we use the convention that repeated
  indices represent summations with respect to these indices,
  and $\partial_j=\partial/\partial x_j$, $j=1,2,3$. Obviously,
$$
  {\mathcal A}(\rho)\in L(h^{m+\theta}(\overline{\Bbb B}^3),
  h^{m-2+\theta}(\overline{\Bbb B}^3)),
  \qquad\vec{\mathcal B}(\rho)\in L(h^{m+\theta}(\overline
  {\Bbb B}^3),(h^{m-1+\theta}(\overline{\Bbb B}^3))^3),
$$
$$
  \vec{\mathcal B}(\rho)\cdot\in L((h^{m+\theta}(\overline
  {\Bbb B}^3))^3, h^{m-1+\theta}(\overline{\Bbb B}^3)),
  \qquad
  \vec{\mathcal B}(\rho)\otimes\in L\big((h^{m+\theta}(\overline
  {\Bbb B}^3))^3,(h^{m-1+\theta}(\overline{\Bbb B}^3))^{3\times
  3}\big).
$$
  The definitions (2.4)--(2.7) can be respectively briefly rewritten
  as follows:
$$
  {\mathcal A}(\rho)u=(\Delta(u\circ\Phi_\rho^{-1}))\circ\Phi_\rho,
  \quad
  \vec{\mathcal B}(\rho)u=(\nabla(u\circ\Phi_\rho^{-1}))\circ\Phi_\rho,
$$
$$
  \vec{\mathcal B}(\rho)\cdot{\bfv}=(\nabla\cdot({\bfv}\circ\Phi_\rho^{-1}))
  \circ\Phi_\rho, \quad
  \vec{\mathcal B}(\rho)\otimes{\bfv}=(\nabla\otimes({\bfv}\circ\Phi_\rho^{-1}))
  \circ\Phi_\rho.
$$
  As in \cite{Fried1} we introduce the following vector functions:
$$
  \bfw_1(x)=(0,x_3,-x_2),\quad \bfw_2(x)=(-x_3,0,x_1),\quad
  \bfw_3(x)=(x_2,-x_1,0).
$$
  Then clearly $\bfv\times x=(\bfv\cdot\bfw_1,\bfv\cdot\bfw_2,\bfv\cdot\bfw_3)$.

  Let $\bfn$ and $\kappa$ be respectively the unit outward normal and the
  mean curvature of $\Gamma_\rho$ (see (1.5)). We denote
$$
  \widetilde{\bfn}_\rho(x)=\bfn(\phi_\rho(x)) \quad\mbox{and}\quad
  \widetilde{\kappa}_\rho(x)=\kappa(\phi_\rho(x)), \quad
  \mbox{for}\;\; x\in {\Bbb S}^2.
$$
  A direct computation shows that
$$
  \widetilde{\bfn}_\rho(x)=\frac{x\cdot[(D\Phi_\rho(x))^{-1}\big]^T}
  {|x\cdot[(D\Phi_\rho(x))^{-1}\big]^T|}=
  {a_{ij}^\rho(x)x_j\bfe_i\over\big|a_{ij}^\rho(x)x_j\bfe_i\big|}
  \equiv(\widetilde{\bfn}_\rho^1(x),\widetilde{\bfn}_\rho^2(x),
  \widetilde{\bfn}_\rho^3(x)),
\eqno{(2.8)}
$$
  where
$$
  \bfe_1=(1,0,0),\quad \bfe_2=(0,1,0),\quad \bfe_3=(0,0,1),
$$
  and
$$
  \widetilde{\kappa}_\rho(x)={1\over2}a_{ij}^\rho(x)\partial_j\widetilde
  {\bfn}_\rho^i(x).
\eqno{(2.9)}
$$

  We remind the reader to notice that, obviously, $[\rho\to\Phi_\rho]
  \in C^\infty(O^{m+\theta}_\delta({\Bbb S}^2),{\rm Diff}^{m+\theta}_h
  ({\Bbb R}^3,{\Bbb R}^3))$. Thus we have
$$
\left\{
\begin{array}{l}
  [\rho\to a_{ij}^\rho]\in C^\infty(O^{m+\theta}_\delta({\Bbb S}^2),
  h^{m-1+\theta}({\Bbb R}^3)),\quad i,j=1,2,3,\\ [0.2cm]
  [\rho\to G_\rho]\in C^\infty(O^{m+\theta}_\delta({\Bbb S}^2),
  h^{m-1+\theta}({\Bbb R}^3)),\\ [0.2cm]
  [\rho\to H_\rho]\in C^\infty(O^{m+\theta}_\delta({\Bbb S}^2),
  h^{m-1+\theta}({\Bbb S}^2)),\\ [0.2cm]
  [\rho\to\widetilde{\kappa}_\rho]\in C^\infty(O^{m+\theta}_\delta({\Bbb S}^2),
  h^{m-2+\theta}({\Bbb S}^2)),\\ [0.2cm]
  [\rho\to\widetilde{\bfn}_\rho]\in C^\infty(O^{m+\theta}_\delta({\Bbb S}^2),
  (h^{m-1+\theta}({\Bbb S}^2))^3).
\end{array}
\right.
\eqno{(2.10)}
$$

  Finally, for $\sigma$, $\bfv$ and $p$ as in (1.1)--(1.8), we denote
$$
  \widetilde{\sigma}=\sigma\circ\Phi_\rho, \quad
  \widetilde{\bfv}=\bfv\circ\Phi_\rho, \quad
  \widetilde{p}=p\circ\Phi_\rho.
$$
  We also denote $\widetilde{\bfw}_j^\rho=\bfw_j\circ\Phi_\rho$, $j=1,2,3$.

  Using these notations, we see easily that the Hanzawa transformation
  transforms the equations (1.1)--(1.5), (1.7) and (1.8) into the following
  equations, respectively:
$$
  {\mathcal A}(\rho)\widetilde{\sigma}=f(\widetilde{\sigma})
  \quad
  \mbox{in}\;\; \Bbb B^3,\quad t>0,
\eqno{(2.11)}
$$
$$
  \vec{\mathcal B}(\rho)\cdot\widetilde{\bfv}=g(\widetilde{\sigma})
  \quad
  \mbox{in}\;\; \Bbb B^3,\quad t>0,
\eqno{(2.12)}
$$
$$
  -{\mathcal A}(\rho)\widetilde{\bfv}+\vec{\mathcal B}(\rho)
  \widetilde{p}-{1\over 3}\vec{\mathcal B}(\rho)(\vec{\mathcal B}
  (\rho)\cdot\widetilde{\bfv})=0 \quad \mbox{in}\;\; \Bbb B^3,
  \quad t>0,
\eqno{(2.13)}
$$
$$
  \widetilde{\sigma}=1 \quad \mbox{on}\;\; \Bbb S^2,\quad t>0,
\eqno{(2.14)}
$$
$$
  \widetilde{\bfT}_\rho(\widetilde{\bfv},\widetilde{p})
  \widetilde{\bfn}_\rho=-\gamma\widetilde{\kappa}_\rho
  \widetilde{\bfn}_\rho \quad \mbox{on}\;\; \Bbb S^2,
  \quad t>0,
\eqno{(2.15)}
$$
$$
  \int_{|x|<1}\widetilde{\bfv}(x)G_\rho(x)dx=0,\quad t>0,
\eqno{(2.16)}
$$
$$
  \int_{|x|<1}\widetilde{\bfv}(x)\cdot\widetilde{\bfw}_j^\rho(x)
  G_\rho(x)dx=0,\quad j=1,2,3, \quad t>0.
\eqno{(2.17)}
$$
  Here $\widetilde{\bfT}_\rho(\widetilde{\bfv},\widetilde{p})
  =[\vec{\mathcal B}(\rho)\otimes\widetilde{\bfv}+(\vec{\mathcal B}
  (\rho)\otimes\widetilde{\bfv})^T]-[\widetilde{p}+(2/3)
  \vec{\mathcal B}(\rho)\cdot\widetilde{\bfv}]\bfI$.

  In what follows we rewrite (1.6) into an explicit equation expressed with
  the function $\rho=\rho(\omega,t)$. Let $\psi_\rho(x,t)=r-1-\rho(\omega,t)$,
  where $r=|x|$ and $\omega=x/|x|$. Then $x\in\Gamma_{\rho}(t)$ if and only if
  $\psi_\rho(x,t)=0$. It follows that the normal velocity of $\Gamma_\rho(t)$
  is as follows (see \cite{EscSim}):
$$
  V_n(x,t)=\frac{\partial_t\rho(\omega,t)}{|\nabla_x\psi_\rho(x,t)|} \quad
  \mbox{for}\;\; x\in\Gamma_\rho(t), \quad t>0.
$$
  Moreover, $\bfn(x,t)=\nabla_x\psi_\rho(x,t)/|\nabla_x\psi_\rho(x,t)|$.
  Hence (1.6) can be rewritten as follows:
$$
  \partial_t\rho(\omega,t)=\bfv(x,t)\cdot\nabla_x\psi_\rho(x,t) \quad
  \mbox{for}\;\; x\in\Gamma_\rho(t), \quad t>0,
$$
  where $\omega=x/|x|$. Since $\nabla_x\psi_\rho=\displaystyle
  \frac{\partial\psi_\rho}{\partial r}\omega+r^{-1}\nabla_\omega
  \psi_\rho=\omega-r^{-1}\nabla_\omega\rho$, we see that after the
  Hanzawa transformation, this equation has the following form:
$$
  \partial_t\rho(\omega,t)=\widetilde{\bfv}(\omega,t)\cdot\big[
  \omega-{\nabla_\omega\rho(\omega,t)\over1+\rho(\omega,t)}\big]
  \quad
  \mbox{for}\;\; \omega\in\Bbb S^2, \quad t>0.
\eqno{(2.18)}
$$
  Finally, we rewrite (1.13) as follows:
$$
  \rho(\omega,0)=\rho_0(\omega) \quad
  \mbox{for}\;\; \omega\in\Bbb S^2.
\eqno{(2.19)}
$$

  In summary, we have the following preliminary result:
\medskip

  {\bf Lemma 2.1}\ \ {\em If $(\sigma,\bfv,p,\rho)$ is a solution
  of the problem $(1.1)$--$(1.9)$, then by letting $\widetilde{\sigma}
  =\sigma\circ\Phi_\rho$, $\widetilde{\bfv}=\bfv\circ\Phi_\rho$ and
  $\widetilde{p}=p\circ\Phi_\rho$, we have that $(\widetilde{\sigma},
  \widetilde{\bfv},\widetilde{p},\rho)$ is a solution of the problem
  $(2.11)$--$(2.19)$. Conversely, If $(\widetilde{\sigma},\widetilde{\bfv},
  \widetilde{p},\rho)$ is a solution of the problem $(2.11)$--$(2.19)$,
  then by letting $\sigma=\widetilde{\sigma}\circ\Phi_\rho^{-1}$,
  $\bfv=\widetilde{\bfv}\circ\Phi_\rho^{-1}$ and $p=\widetilde{p}
  \circ\Phi_\rho^{-1}$, we have that  $(\sigma,\bfv,p,\rho)$ is a
  solution of the problem $(1.1)$--$(1.9)$.}$\quad\Box$
\medskip

  In the sequel we further reduce the problem (2.11)--(2.19) into a scalar
  equation containing the single unknown $\rho$. The idea is to first solve
  the system of equations (2.11)--(2.17) to get $\widetilde{\sigma}$,
  $\widetilde{\bfv}$ and $\widetilde{p}$ as functionals of $\rho$, and next
  substitute $\widetilde{\bfv}$ obtained in this way into the equation (2.18).

  We first consider (2.11) and (2.14). We have:
\medskip

  {\bf Lemma 2.2}\ \ {\em Let $\delta$ be sufficiently small. Then, given $\rho
  \in O^{m+\theta}_\delta(\Bbb S^2)$, the boundary value problem
$$
  {\mathcal A}(\rho)\widetilde{\sigma}=f(\widetilde{\sigma}) \quad
  \mbox{in}\;\; \Bbb B^3, \qquad
  \widetilde{\sigma}=1 \quad \mbox{on}\;\; \Bbb S^2
\eqno{(2.20)}
$$
  has a unique solution $\widetilde{\sigma}={\mathcal R}(\rho)\in h^{m+\theta}
  (\overline{\Bbb B}^3)$ which satisfies $0<\widetilde{\sigma}\leq 1$.
  Moreover, we have
$$
  {\mathcal R}\in C^\infty(O^{m+\theta}_\delta(\Bbb S^2),
  h^{m+\theta}(\overline{\Bbb B}^3)).
\eqno{(2.21)}
$$
}

  {\em Proof}:\ \ See Lemma 3.1 of \cite{CuiEsc1}.$\qquad$ $\Box$
\medskip

  Next, for given $\rho\in O^{m+\theta}_\delta(\Bbb S^2)$ we consider the
  following boundary value problem:
$$
  \vec{\mathcal B}(\rho)\cdot\widetilde{\bfv}=\varphi \quad
  \mbox{in}\;\; \Bbb B^3,
\eqno{(2.22)}
$$
$$
  -{\mathcal A}(\rho)\widetilde{\bfv}+\vec{\mathcal B}(\rho)\widetilde{p}
  =\bfg \quad \mbox{in}\;\; \Bbb B^3,
\eqno{(2.23)}
$$
$$
  \widetilde{\bfT}_\rho(\widetilde{\bfv},\widetilde{p})\widetilde{\bfn}_\rho=
  \bfh \quad \mbox{on}\;\; \Bbb S^2,
\eqno{(2.24)}
$$
$$
  \int_{|x|<1}\widetilde{\bfv}(x)G_\rho(x)dx=0,
\eqno{(2.25)}
$$
$$
  \int_{|x|<1}\widetilde{\bfv}(x)\cdot\widetilde{\bfw}_j^\rho(x)
  G_\rho(x)dx=0,\quad j=1,2,3,
\eqno{(2.26)}
$$
  where $\varphi\in h^{m\!-\!k\!-\!1\!+\!\theta}(\overline{\Bbb B}^3)$, $\bfg
  \in (h^{m\!-\!k\!-\!2\!+\!\theta}(\overline{\Bbb B}^3))^3$ and $\bfh\in
  (h^{m\!-\!k\!-\!1\!+\!\theta}({\Bbb S}^2))^3$ for some $0\!\leq k\leq\!
  m\!-\!2$.
\medskip

  {\bf Lemma 2.3}\ \ {\em Let $\delta$ be sufficiently small and let $\rho\in
  O^{m+\theta}_\delta(\Bbb S^2)$ be given. A necessary and sufficient condition
  for $(2.22)$--$(2.26)$ to have a solution is that $\varphi$, $\bfg$ and $\bfh$
  satisfy the following relations:
$$
  \int_{|x|<1}\big(\bfg(x)-{1\over 3}\vec{\mathcal B}(\rho)\varphi(x)\big)
  \cdot\widetilde{\bfw}_j^\rho(x) G_\rho(x)dx+
  \int_{|x|=1}\bfh(x)\cdot\widetilde{\bfw}_j^\rho(x)
  H_\rho(x) dS_x=0, \quad j=1,2,3,
\eqno{(2.27)}
$$
$$
  \int_{|x|<1}\big(\bfg(x)-{1\over 3}\vec{\mathcal B}(\rho)\varphi(x)\big)
  \cdot{\bfe}_j G_\rho(x)dx+
  \int_{|x|=1}\bfh(x)\cdot{\bfe}_j H_\rho(x) dS_x=0, \quad j=1,2,3.
\eqno{(2.28)}
$$
  If this condition is satisfied, then $(2.12)$--$(2.26)$ has a unique solution
  $(\widetilde{\bfv},\widetilde{p})\in (h^{m-k+\theta}(\overline{\Bbb B}^3))^3\times
  h^{m-k-1+\theta}(\overline{\Bbb B}^3)$.}
\medskip

  {\em Proof}:\ \ Integrating by parts and employing the divergence theorem
  we see that for any $\bfv$, $\bfw\in(C^2(\overline{\Omega}_\rho))^3$ and
  $p\in C^1(\overline{\Omega}_\rho))$ there holds the following integral
  identity (cf. \cite{GunProk}):
$$
\begin{array}{rcl}
  &&\displaystyle{1\over 2}\int_{{\Omega}_\rho}[\nabla\otimes\bfv+
  (\nabla\otimes\bfv)^T]\cdot[\nabla\otimes\bfw+(\nabla\otimes\bfw)^T]dx
  -\int_{{\Omega}_\rho}[p+{2\over 3}(\nabla\cdot\bfv)]
  \nabla\cdot\bfw dx  \\ [0.2cm]
  &=&\displaystyle\int_{{\Omega}_\rho}[-\Delta\bfv+\nabla p-
  {1\over 3}\nabla(\nabla\cdot\bfv)]\cdot\bfw dx
  +\int_{\Gamma_\rho}\bfT(\bfv,p)\bfn\cdot\bfw dS_x.
\end{array}
\eqno{(2.29)}
$$
  Here, for two matrix $A$ and $B$ we use the notation $A\cdot B$ to denote the
  inner product of $A$ and $B$, i.e., $A\cdot B=A_{ij}B_{ij}$. By
  making the Hanzawa transformation, we see that for any $\widetilde{\bfv}$,
  $\widetilde{\bfw}\in (C^2(\overline{\Bbb B}^3))^3$ and $\widetilde{p}\in
  C^1(\overline{\Bbb B}^3)$ there holds
$$
\begin{array}{rcl}
  &&\displaystyle{1\over 2}\!\int_{|x|\!<\!1}\!\!
  [\vec{\mathcal B}(\rho)\!\otimes\!\widetilde{\bfv}\!+\!
  (\vec{\mathcal B}(\rho)\!\otimes\!\widetilde{\bfv})^T]\!\cdot\!
  [\vec{\mathcal B}(\rho)\!\otimes\!\widetilde{\bfw}\!+\!
  (\vec{\mathcal B}(\rho)\!\otimes\!\widetilde{\bfw})^T]G_\rho dx
%\\ [0.2cm]&&\displaystyle
  \!-\!\int_{|x|\!<\!1}\!\![\widetilde{p}+
  {2\over 3}\vec{\mathcal B}(\rho)\!\cdot\!\widetilde{\bfv}]
  \vec{\mathcal B}(\rho)\!\cdot\!\widetilde{\bfw}G_\rho dx
\\ [0.3cm]
  &=&\displaystyle\int_{|x|<1}[-{\mathcal A}(\rho)\widetilde{\bfv}+
  \vec{\mathcal B}(\rho)\widetilde{p}-{1\over 3}
  \vec{\mathcal B}(\rho)(\vec{\mathcal B}(\rho)\cdot\widetilde{\bfv})]
  \cdot\widetilde{\bfw}G_\rho dx
  +\int_{|x|=1}\widetilde{\bfT}_\rho(\widetilde{\bfv},\widetilde{p})
  \bfn_\rho\cdot\widetilde{\bfw} H_\rho dS_x.
\end{array}
\eqno{(2.30)}
$$
  Besides, clearly $\nabla\otimes\bfw_j+(\nabla\otimes\bfw_j)^T=0$,
  $\nabla\cdot\bfw_j=0$, $j=1,2,3$, which yield, after the Hanzawa
  transformation, that
$$
  \vec{\mathcal B}(\rho)\otimes\widetilde{\bfw}_j^\rho+
  (\vec{\mathcal B}(\rho)\otimes\widetilde{\bfw}_j^\rho)^T=0, \quad
  \vec{\mathcal B}(\rho)\cdot\widetilde{\bfw}_j^\rho=0, \quad
  j=1,2,3.
$$
  Hence, if $(\widetilde{\bfv},\widetilde{p})$ is a solution of (2.22)--(2.26),
  then by takeing $\widetilde{\bfw}=\widetilde{\bfw}_j^\rho$ ($j=1,2,3$) in
  (2.30), we see that (2.27) holds. Similarly we have (2.28). This proves the
  necessity of (2.27) and (2.28).

  Next we assume that the conditions (2.27) and (2.28) are satisfied,
  and proceed to prove that there exists a unique solution to the
  problem (2.22)--(2.26). We first prove uniqueness of the solution.
  Let $(\widetilde{\bfv}_1,\widetilde{p}_1)$ and $(\widetilde{\bfv}_2,
  \widetilde{p}_2)$ be two solutions of (2.22)--(2.26). Then
  $\widetilde{\bfv}=\widetilde{\bfv}_1-\widetilde{\bfv}_2$ and
  $\widetilde{p}=\widetilde{p}_1-\widetilde{p}_2$ satisfy the
  corresponding homogeneous equations. Thus, by letting
  $\widetilde{\bfw}=\widetilde{\bfv}$ in (2.30), we get
$$
  \int_{|x|<1}[\vec{\mathcal B}(\rho)\otimes\widetilde{\bfv}+
  (\vec{\mathcal B}(\rho)\otimes\widetilde{\bfv})^T]\cdot
  [\vec{\mathcal B}(\rho)\otimes\widetilde{\bfv}+
  (\vec{\mathcal B}(\rho)\otimes\widetilde{\bfv})^T]G_\rho(x)dx=0,
$$
  so that
$$
  \vec{\mathcal B}(\rho)\otimes\widetilde{\bfv}+
  (\vec{\mathcal B}(\rho)\otimes\widetilde{\bfv})^T=0 \quad
  \mbox{in}\;\; \Bbb B^3.
$$
  This combined with (2.25) and (2.26) yields, by the Korn inequality
$$
  \|\bfu\|^2_{(H^1(\Omega))^3}\leq C_1\|{\mathbf S}(\bfu)\|^2_{(L^2(\Omega))^{3\times3}}
  +C_2\left(\left|\int_\Omega\bfu dx\right|^2+
  \left|\int_\Omega\bfu\times x dx\right|^2\right),
\eqno{(2.31)}
$$
  where ${\mathbf S}(\bfu)=\nabla\otimes\bfu+(\nabla\otimes\bfu)^T$
  (cf. Proposition 8.1 of \cite{Solon}), that $\widetilde{\bfv}=0$.
  From this it follows immediately that also $\widetilde{p}=0$. Hence
  the solution is unique if it exists.

  To prove existence we denote, for a given $0\le k\le m-2$,
$$
  {\mathbb X}=(h^{m-k+\theta}(\overline{\Bbb B}^3))^3\times
  h^{m-k-1+\theta}(\overline{\Bbb B}^3)\times {\Bbb R}^6,
$$
$$
  {\mathbb Y}=h^{m-k-1+\theta}(\overline{\Bbb B}^3)\times
  (h^{m-k-2+\theta}(\overline{\Bbb B}^3))^3\times
  (h^{m-k-1+\theta}(\Bbb S^2))^3\times {\Bbb R}^3\times {\Bbb R}^3,
$$
  and regard $O^{m+\theta}_\delta(\Bbb S^2)$ as an open subset of the Banach
  space $h^{m+\theta}(\Bbb S^2)$. For every $\rho\in O^{m+\theta}_\delta
  (\Bbb S^2)$, we define a linear operator ${\mathcal L}(\rho): {\mathbb X}\to
  {\mathbb Y}$ as follows:
$$
  {\mathcal L}(\rho)U=\left[
\begin{array}{c}
  \vec{\mathcal B}(\rho)\cdot\widetilde{\bfv}\\ [0.2cm]
  -{\mathcal A}(\rho)\widetilde{\bfv}+\vec{\mathcal B}
  (\rho)\widetilde{p}+l_\rho(\zeta)\\ [0.2cm]
  \widetilde{\bfT}_\rho(\widetilde{\bfv},\widetilde{p})
  \widetilde{\bfn}_\rho\\ [0.2cm]
  \displaystyle\int_{|x|<1}\widetilde{\bfv}(x)G_\rho(x)dx
  \\ [0.3cm]
  \displaystyle\int_{|x|<1}\widetilde{\bfv}(x)\cdot
  \widetilde{\bfw}_j^\rho(x)G_\rho(x)dx
\end{array}
\right]^T
  \quad \mbox{for}\;\; U=(\widetilde{\bfv},\widetilde{p},\zeta)
  \in {\mathbb X},
$$
  where $l_\rho$ is the linear operator from $\Bbb R^6$ to $(h^{m-k-2+\theta}
  (\overline{\Bbb B}^3))^3$ defined by $l(\zeta)=\bfa+b_1\widetilde{\bfw}_1^\rho+
  b_2\widetilde{\bfw}_2^\rho+b_3\widetilde{\bfw}_3^\rho$ for $\zeta=(\bfa,b_1,
  b_2,b_3)\in\Bbb R^3\times\Bbb R\times\Bbb R\times\Bbb R$. By (2.10) it is clear
  that
$$
  {\mathcal L}\in C^\infty(O^{m+\theta}_\delta(\Bbb S^2),
  L({\mathbb X},{\mathbb Y})),
$$
  and
$$
  {\mathcal L}(0)U=\left[
\begin{array}{c}
  \nabla\cdot{\bfv}\\ [0.2cm]
  -\Delta\bfv+\nabla {p}+l_0(\zeta)\\ [0.2cm]
  \bfT(\bfv,p){\bfn}\\ [0.2cm]
  \displaystyle\int_{|x|<1}\bfv(x)dx\\ [0.3cm]
  \displaystyle\int_{|x|<1}\bfv(x)\times x dx
\end{array}
\right]^T
\quad \mbox{for}\;\; U=(\bfv,p,\zeta)\in {\mathbb X},
$$
  where $l_0$ is the linear operator from $\Bbb R^6$ to $(h^{m-k-2+\theta}
  (\overline{\Bbb B}^3))^3$ defined by $l(\zeta)=\bfa+b_1\bfw_1+b_2\bfw_2+b_3\bfw_3$ for
  $\zeta=(\bfa,b_1,b_2,b_3)\in\Bbb R^3\times\Bbb R\times\Bbb R\times\Bbb R$.
  From the proof of Theorem 2.1 of \cite{Fried1} we see that ${\mathcal L}(0)$
  is an isomorphism from ${\mathbb X}$ to ${\mathbb Y}$ (cf. also Lemma A.1 of
  \cite{EscProk}). Since all isomorphisms from ${\mathbb X}$ to ${\mathbb Y}$
  forms an open set in $L({\mathbb X},{\mathbb Y})$, we conclude that for
  $\delta$ sufficiently small, ${\mathcal L}(\rho)$ is also an isomorphism from
  ${\mathbb X}$ to ${\mathbb Y}$ for any $\rho\in O^{m+\theta}_\delta(\Bbb S^2)$.
  This particularly implies that given $\varphi\in h^{m-k-1+\theta}(\overline
  {\Bbb B}^3)$, $\bfg\in (h^{m-k-2+\theta}(\overline{\Bbb B}^3))^3$ and $\bfh\in
  (h^{m-k-1+\theta}(\overline{\Bbb B}^3))^3$, there exist unique $\widetilde{\bfv}
  \in(h^{m-k+\theta}(\overline{\Bbb B}^3))^3$, $\widetilde{p}\in h^{m-k-1+\theta}
  (\overline{\Bbb B}^3)$ and $\zeta\in {\Bbb R}^6$ such that they satisfy
  (2.22), (2.24)--(2.26) and
$$
  -{\mathcal A}(\rho)\widetilde{\bfv}+\vec{\mathcal B}(\rho)\widetilde{p}
  +l_\rho(\zeta)=\bfg \quad \mbox{in}\;\; \Bbb B^3.
\eqno{(2.32)}
$$
  We claim that $\zeta=0$. Indeed, taking $\widetilde{\bfw}=\widetilde{\bfw}_j,
  \widetilde{\bfe}_j$ in (2.30) and using (2.22), (2.24), (2.27), (2.28) and
  (2.32), we get
$$
  \int_{|x|<1}l_\rho(\zeta)\cdot\widetilde{\bfw}_j^\rho G_\rho dx=0, \quad
  \int_{|x|<1}l_\rho(\zeta)\cdot{\bfe}_j G_\rho dx=0, \quad j=1,2,3.
$$
  From these relations we can easily show that if $\rho=0$ then $\zeta=0$.
  By continuity (a small perturbation of a nonsingular matrix is still
  nonsingular), this implies that if $\delta$ is sufficiently small then
  for any $\rho\in O^{m+\theta}_\delta({\Bbb S}^2)$ we also have $\zeta=0$.
  Hence our claim is true. It follows that $(\widetilde{\bfv},\widetilde{p})$
  is a solution of (2.22)--(2.26). This completes the proof of Lemma 2.3.
  $\qquad$ $\Box$
\medskip

  {\bf Lemma 2.4}\ \ {\em For the solution of the problem $(2.22)$--$(2.26)$,
  we have $\widetilde{\bfv}=\vec{\mathcal P}(\rho)\varphi+{\bf Q}(\rho)\bfg
  +{\bf R}(\rho)\bfh$, where
$$
\left\{
\begin{array}{l}
\displaystyle
  \vec{\mathcal P}\in\bigcap_{k=0}^{m-2} C^\infty(O^{m+\theta}_\delta(\Bbb S^2),
  L(h^{m-k-1+\theta}(\overline{\Bbb B}^3),
  (h^{m-k+\theta}(\overline{\Bbb B}^3))^3),\\ [0.3cm]
\displaystyle
  {\bf Q}\in\bigcap_{k=0}^{m-2} C^\infty(O^{m+\theta}_\delta(\Bbb S^2),
  L((h^{m-k-2+\theta}(\overline{\Bbb B}^3))^3,
  (h^{m-k+\theta}(\overline{\Bbb B}^3))^3)),\\ [0.3cm]
\displaystyle
  {\bf R}\in\bigcap_{k=0}^{m-2} C^\infty(O^{m+\theta}_\delta(\Bbb S^2),
  L((h^{m-k-1+\theta}({\Bbb S}^2))^3,
  (h^{m-k+\theta}(\overline{\Bbb B}^3))^3)).
\end{array}
\right.
$$
}

  {\em Proof}:\ \ Let notations be as in the proof of Lemma 2.3. We denote by
  ${\mathcal I}_1$, ${\mathcal I}_2$ and ${\mathcal I}_3$ the natural embedding
  operators from $h^{m-k-1+\theta}(\overline{\Bbb B}^3)$, $(h^{m-k-2+\theta}
  (\overline{\Bbb B}^3))^3$ and $(h^{m-k-1+\theta}({\Bbb S}^2))^3$, respectively,
  into ${\mathbb Y}$, and by ${\mathcal J}$ the projection operator from
  ${\mathbb X}$ onto $(h^{m-k+\theta}(\overline{\Bbb B}^3))^3$. Then by letting
$$
  \vec{\mathcal P}(\rho)={\mathcal J}\circ {\mathcal L}(\rho)^{-1}
  \circ{\mathcal I}_1, \quad
  {\bf Q}(\rho)={\mathcal J}\circ {\mathcal L}(\rho)^{-1}\circ{\mathcal I}_2,
  \quad
  {\bf R}(\rho)={\mathcal J}\circ {\mathcal L}(\rho)^{-1}\circ{\mathcal I}_3,
$$
  we immediately see that the desired assertion follows. $\qquad$ $\Box$
\medskip

  The system of equations (2.12), (2.13) and (2.15)--(2.17) can be rewritten in
  the form of (2.22)--(2.26), with
$$
  \varphi=g(\widetilde{\sigma}), \quad
  \bfg={1\over 3}\vec{\mathcal B}(\rho)g(\widetilde{\sigma}), \quad
  \bfh=-\gamma\widetilde{\kappa}_\rho\widetilde{\bfn}_\rho.
\eqno{(2.33)}
$$
  We assert that (2.27) and (2.28) are satisfied by these functions. Indeed,
  since $\bfg={1\over 3}\vec{\mathcal B}(\rho)\varphi$, this assertion follows
  if we show that
$$
  \int_{|x|=1}\bfh(x)\cdot\widetilde{\bfw}_j^\rho(x)
  H_\rho(x) dS_x=0, \quad
  \int_{|x|=1}\bfh(x)\cdot{\bfe}_j H_\rho(x) dS_x=0, \quad j=1,2,3.
\eqno{(2.34)}
$$
  Let $\Delta_{\Gamma_\rho}$ be the Laplace-Beltrami operator on $\Gamma_\rho$.
  Then we have $\kappa(x)\bfn(x)=-\Delta_{\Gamma_\rho}x$ for $x\in\Gamma_\rho$
  (cf. \cite{EscProk}, \cite{Solon}). Since $\Delta_{\Gamma_\rho}$ is a
  symmetric operator in $L^2(\Gamma_\rho,dS_x)$, we see that
$$
  \int_{\Gamma_\rho}(\Delta_{\Gamma_\rho}x_i\cdot x_j-
  \Delta_{\Gamma_\rho}x_j\cdot x_i)dS_x=0, \quad i,j=1,2,3.
\eqno{(2.35)}
$$
  Thus
$$
   \int_{|x|=1}\widetilde{\kappa}_\rho\bfn_\rho\cdot\widetilde{\bfw}_j^\rho
  H_\rho dS_x=\int_{\Gamma_\rho}\kappa\bfn\cdot\bfw_j dS_x
  =-\int_{\Gamma_\rho}\Delta_{\Gamma_\rho}x\cdot\bfw_j dS_x=0, \quad j=1,2,3.
$$
  Similarly we have
  $$
   \int_{|x|=1}\widetilde{\kappa}_\rho\bfn_\rho\cdot{\bfe}_j H_\rho dS_x
   =\int_{\Gamma_\rho}\kappa\bfn\cdot{\bfe}_j dS_x
  =-\int_{\Gamma_\rho}\Delta_{\Gamma_\rho}x\cdot{\bfe}_j dS_x=0, \quad j=1,2,3.
$$
  This verifies (2.34).

  Now, given $\rho\in O^{m+\theta}_\delta(\Bbb S^2)$, we first use Lemma 2.2
  to solve the equations (2.11) and (2.14). This gives $\widetilde{\sigma}=
  {\mathcal R}(\rho)\in h^{m+\theta}(\overline{\Bbb B}^3)$. Next we use Lemma
  2.3 to solve equations (2.12), (2.13) and (2.15)--(2.17). Note that with
  $\varphi$, $\bfg$ and $\bfh$ given in (2.33), we have $\varphi\in h^{m+\theta}
  (\overline{\Bbb B}^3)\subseteq h^{m-2+\theta}(\overline{\Bbb B}^3)$,
  $\bfg\in (h^{m-1+\theta}(\overline{\Bbb B}^3))^3\subseteq (h^{m-3+\theta}
  (\overline{\Bbb B}^3))^3$, and, by (2.10), $\bfh\in (h^{m-2+\theta}
  ({\Bbb S}^2))^3$. Hence, by Lemma 2.3 (with $k=1$) it follows that these
  equations have a unique solution $(\widetilde{\bfv},\widetilde{p})\in
  (h^{m-1+\theta}(\overline{\Bbb B}^3))^3\times h^{m-2+\theta}
  (\overline{\Bbb B}^3)$. Moreover, since $\widetilde{\sigma}=
  {\mathcal R}(\rho)$, by Lemma 2.4 we have
$$
  \widetilde{\bfv}=\vec{\mathcal P}(\rho)g({\mathcal R}(\rho))
  +{1\over 3}{\mathbf Q}(\rho)\vec{\mathcal B}(\rho)g({\mathcal R}(\rho))
  -\gamma {\mathbf R}(\rho)({\mathcal K}(\rho)\vec{\mathcal N}(\rho)).
\eqno{(2.36)}
$$
  where ${\mathcal K}(\rho)=\widetilde{\kappa}_\rho$ and $\vec{\mathcal N}
  (\rho)=\widetilde{\bfn}_\rho$. We note that
$$
  {\mathcal K}\in C^\infty(O^{m+\theta}_\delta({\Bbb S}^2),
  h^{m-2+\theta}({\Bbb S}^2)), \quad
  \vec{\mathcal N}\in C^\infty(O^{m+\theta}_\delta({\Bbb S}^2),
  (h^{m-1+\theta}({\Bbb S}^2))^3).
\eqno{(2.37)}
$$
  Substituting the expression of $\widetilde{\bfv}$ in (2.36) into (2.18),
  and introducing the operator ${\mathcal Q}:O^{m+\theta}_\delta({\Bbb S}^2)
  \to h^{m-1+\theta}({\Bbb S}^2)$ by
$$
  {\mathcal Q}(\rho)={\rm tr}_{{\Bbb S}^2}\big[\vec{\mathcal P}(\rho)
  g({\mathcal R}(\rho))+{1\over 3}{\mathbf Q}(\rho)\vec{\mathcal B}(\rho)
  g({\mathcal R}(\rho))-\gamma {\mathbf R}(\rho)({\mathcal K}(\rho)
  \vec{\mathcal N}(\rho))\big]\cdot\big[\omega-{\nabla_\omega\rho
  \over1+\rho}\big],
\eqno{(2.38)}
$$
  (for $\rho\in O^{m+\theta}_\delta({\Bbb S}^2)$), where as before
  $\omega$ represents the variable in ${\Bbb S}^2$, we see that the
  problem (2.11)--(2.19) is reduced into the following initial value
  problem for a differential equation in the Banach space
  $h^{m-1+\theta}({\Bbb S}^2)$:
$$
\left\{
\begin{array}{l}
   \partial_t\rho={\mathcal Q}(\rho), \quad t>0,\\ [0.1cm]
   \rho|_{t=0}=\rho_0.
\end{array}
\right.
\eqno{(2.39)}
$$
  Here ${\mathcal Q}$ is regarded as a unbounded operator in $h^{m-1+\theta}
  ({\Bbb S}^2)$ with domain $O^{m+\theta}_\delta({\Bbb S}^2)$.

  We summarize:
\medskip

  {\bf Lemma 2.5}\ \ {\em Let $(\widetilde{\sigma},\widetilde{\bfv},
  \widetilde{p},\rho)$ be a solution of the problem $(2.11)$--$(2.19)$.
  Then $\rho$ is a solution of the initial value problem $(2.39)$.
  Conversely, if $\rho$ is a solution of the initial value problem
  $(2.39)$, then by letting $\widetilde{\sigma}={\mathcal R}(\rho)$
  and $(\widetilde{\bfv},\widetilde{p})=P{\mathcal L}(\rho)^{-1}
  (g(\widetilde{\sigma}),{1\over 3}\vec{\mathcal B}(\rho)g(\widetilde
  {\sigma}),-\gamma\widetilde{\kappa}_\rho\widetilde{\bfn}_\rho,0,0)$,
  where $P$ denotes the projection from ${\mathbb X}=(h^{m-1+\theta}
  (\overline{\Bbb B}^3))^3\times h^{m-2+\theta}(\overline{\Bbb B}^3)
  \times {\Bbb R}^6$ onto $(h^{m-1+\theta}(\overline{\Bbb B}^3))^3
  \times h^{m-2+\theta}(\overline{\Bbb B}^3)$, we have that
  $(\widetilde{\sigma},\widetilde{\bfv},\widetilde{p},\rho)$
  is a solution of  $(2.11)$--$(2.19)$.} $\qquad$ $\Box$
\medskip

  From (2.10), (2.21), (2.37) and Lemma 2.4 we see that
$$
  {\mathcal Q}\in C^\infty(O_\delta^{m+\theta}
  ({\Bbb S}^2),h^{m-1+\theta}({\Bbb S}^2)).
\eqno{(2.40)}
$$
  In the sequel we prove that if $\delta$ is sufficiently small then for any
  $\rho\in O_\delta^{m+\theta}({\Bbb S}^2)$, $D{\mathcal Q}(\rho)$ is a
  infinitesimal generator of an analytic semigroup in $h^{m-1+\theta}
  ({\Bbb S}^2)$ with domain $h^{m+\theta}({\Bbb S}^2)$, so that the
  differential equation in (2.39) is of parabolic type. Here and in what
  follows, the notation $D\cdot$ represents Fr\'{e}chet derivatives of smooth
  operators from $h^{m+\theta}({\Bbb S}^2)$ to $h^{m-1+\theta}({\Bbb S}^2)$.

  We first note that the mean curvature operator ${\mathcal K}$ has the
  following expression
$$
  {\mathcal K}(\rho)={\mathcal K}_1(\rho)\rho+{\mathcal K}_0(\rho),
\eqno{(2.41)}
$$
  where, for each $\rho$, ${\mathcal K}_1(\rho)$ is a second-order linear
  elliptic partial differential operator on ${\Bbb S}^2$ with coefficients
  being functions of $\rho$ and its first-order derivatives, and
  ${\mathcal K}_0$ is a first-order nonlinear partial differential operator
  on ${\Bbb S}^2$, so that
$$
  {\mathcal K}_1\in C^\infty(O_\delta^{m+\theta}({\Bbb S}^2),
  L(h^{m+\theta}({\Bbb S}^2),h^{m-2+\theta}({\Bbb S}^2))), \quad
  {\mathcal K}_0\in C^\infty
  (O_\delta^{m+\theta}({\Bbb S}^2),h^{m-1+\theta}({\Bbb S}^2)).
\eqno{(2.42)}
$$
  (see \cite{CuiEsc1} and \cite{EscSim}). Substituting (2.41) into
  (2.38) we see that
$$
  {\mathcal Q}(\rho)={\mathcal Q}_2(\rho)[\rho,\rho]+
  {\mathcal Q}_1(\rho)\rho+{\mathcal Q}_0(\rho),
\eqno{(2.43)}
$$
  where, for each $\rho$, ${\mathcal Q}_2(\rho)$ is a bilinear operator,
  ${\mathcal Q}_1(\rho)$ is a linear operator, and ${\mathcal Q}_0$ is a
  nonlinear operator; they are respectively defined as follows:
$$
\begin{array}{rl}
  &\displaystyle{\mathcal Q}_2(\rho)[\eta_1,\eta_2]=
  \gamma {\rm tr}_{{\Bbb S}^2}\big\{{\mathbf R}(\rho)
  \big[{\mathcal K}_1(\rho)\eta_1\vec{\mathcal N}(\rho)\big]
  \big\}\cdot{\nabla_\omega\eta_2\over1+\rho},\\ [0.3cm]
&\begin{array}{rl}
  {\mathcal Q}_1(\rho)\eta=&\displaystyle-{\rm tr}_{{\Bbb S}^2}\big
  \{\vec{\mathcal P}(\rho)g({\mathcal R}(\rho))
  +{1\over 3}{\mathbf Q}(\rho)\vec{\mathcal B}
  (\rho)g({\mathcal R}(\rho))-\gamma {\mathbf R}(\rho)
  \big[{\mathcal K}_0(\rho)\vec{\mathcal N}(\rho)\big]\big\}
  \cdot{\nabla_\omega\eta\over1+\rho}\\ [0.2cm]
  &\displaystyle
  -\gamma {\rm tr}_{{\Bbb S}^2}\big\{{\mathbf R}(\rho)
  \big[{\mathcal K}_1(\rho)\eta
  \vec{\mathcal N}(\rho))\big]\big\}\cdot\omega,
\end{array}\\ [0.4cm]
  &\displaystyle
  {\mathcal Q}_0(\rho)={\rm tr}_{{\Bbb S}^2}\big\{\vec{\mathcal P}(\rho)
  g({\mathcal R}(\rho))+{1\over 3}{\mathbf Q}(\rho)
  \vec{\mathcal B}(\rho)g({\mathcal R}(\rho))-\gamma {\mathbf R}(\rho)
  \big[{\mathcal K}_0(\rho)\vec{\mathcal N}(\rho)\big]\big\}\cdot\omega.
\end{array}
$$
  We note that, by (2.21), (2.37), (2.42) and Lemma 2.4 we have
$$
\begin{array}{l}
  {\mathcal Q}_2\in C^\infty(O_\delta^{m+\theta}({\Bbb S}^2),
  BL(h^{m+\theta}({\Bbb S}^2)\times h^{m+\theta}({\Bbb S}^2),
  h^{m-1+\theta}({\Bbb S}^2))),\\ [0.2cm]
  {\mathcal Q}_1\in C^\infty(O_\delta^{m+\theta}({\Bbb S}^2),
  L(h^{m+\theta}({\Bbb S}^2),h^{m-1+\theta}({\Bbb S}^2))),\\ [0.2cm]
  {\mathcal Q}_0\in C^\infty(O_\delta^{m+\theta}({\Bbb S}^2),
  h^{m+\theta}({\Bbb S}^2)),
\end{array}
$$
  where $BL(\cdot\times\cdot,\cdot)$ denotes the Banach space of all bilinear
  operators with respect the corresponding spaces.

  Given two Banach spaces $E_0$ and $E_1$ such that $E_1$ is continuously and
  densely embedded into $E_0$, we denote by $\mathcal H(E_1,E_0)$ the subset of
  all linear operators $A\in L(E_1,E_0)$ such that $-A$ generates a strongly
  continuous analytic semigroup on $E_0$.
\medskip

  {\bf Lemma 2.6}\ \ $-D{\mathcal Q}(0)\in\mathcal H(h^{m+\theta}(\Bbb S^2),
  h^{m-1+\theta}(\Bbb S^2))$.
\medskip

  {\em Proof:}\ \ For any $\rho\in O_\delta^{m+\theta}({\Bbb S}^2)$ and
  $\eta\in h^{m+\theta}({\Bbb S}^2)$ we have
$$
  D{\mathcal Q}(\rho)\eta={\mathcal Q}_2(\rho)[\eta,\rho]+
  {\mathcal Q}_2(\rho)[\rho,\eta]+
  [D{\mathcal Q}_2(\rho)\eta][\rho,\rho]
  +{\mathcal Q}_1(\rho)\eta+
  [D{\mathcal Q}_1(\rho)\eta]\rho+D{\mathcal Q}_0(\rho)\eta.
$$
  In particular,
$$
  D{\mathcal Q}(0)\eta={\mathcal Q}_1(0)\eta+D{\mathcal Q}_0(0)
  \eta \quad \mbox{for}\;\;\eta\in h^{m+\theta}({\Bbb S}^2),
\eqno{(2.44)}
$$
  i.e., $D{\mathcal Q}(0)={\mathcal Q}_1(0)+D{\mathcal Q}_0(0)$. We note that
  ${\mathcal Q}_1(0)\in L(h^{m+\theta}({\Bbb S}^2),h^{m-1+\theta}({\Bbb S}^2))$
  and $D{\mathcal Q}_0(0)\in L(h^{m+\theta}({\Bbb S}^2),h^{m+\theta}({\Bbb S}^2))$.
  Thus, by a standard perturbation result for infinitesimal generators of continuous
  analytic semigroups (see \cite{Amann} and \cite{Lunar}), the desired assertion
  follows if we prove that $-{\mathcal Q}_1(0)\in\mathcal H(h^{m+\theta}(\Bbb S^2),
  h^{m-1+\theta}(\Bbb S^2))$.

  Since (cf. \cite{FriedReit2})
$$
  \mathcal K(\epsln\eta)=1-\epsln[\eta(\omega)+{1\over2}\Delta_\omega
  \eta(\omega)]+o(\epsln),
$$
  where $\Delta_\omega$ is the Laplace-Beltrami operator on the
  sphere $\Bbb S^2$, we have ${\mathcal K}_0(0)={\mathcal K}(0)=1$
  and ${\mathcal K}_1(0)\eta=-{1\over2}\Delta_{\omega}\eta$. Hence,
  from the definition of ${\mathcal Q}_1$ we see that
$$
  {\mathcal Q}_1(0)\eta=-\bfv|_{\Bbb S^2}\cdot\nabla_\omega\eta
  -\gamma\bfu_{\eta}|_{\Bbb S^2}\cdot\omega,
\eqno{(2.45)}
$$
  where $\bfv$ is the solution
  of the following boundary value problem:
$$
   \Delta\sigma=f(\sigma) \quad\mbox{in}\;\; |x|<1,
$$
$$
  \nabla\cdot\bfv=g(\sigma) \quad\mbox{in}\;\; |x|<1,
$$
$$
   -\Delta\bfv+\nabla p-{1\over3}\nabla
   (\nabla\cdot\bfv)=0 \quad\mbox{in} \;\;|x|<1,
$$
$$
   \sigma=1 \quad\mbox{on}\;\;|x|=1,
$$
$$
  \bfT({\bfv},p){\bfn}=-\gamma{\bfn}
  \quad\mbox{on}\;\;|x|=1,
$$
$$
   \int_{|x|<1}\bfv\;dx=0,
$$
$$
   \int_{|x|<1}\bfv\times{x}d\;x=0,
$$
  and $\bfu_{\eta}=-{1\over2}{\bf R}(0)(\Delta_{\omega}\eta\cdot{\bfn})$.
  Clearly, $\bfv=\bfv_s$ --- the radially symmetric stationary solution of
  the problem (1.1)--(1.8) (see Appendix A), so that $\bfv|_{\Bbb S^2}=0$.
  It follows that
$$
  {\mathcal Q}_1(0)\eta=-\gamma\bfu_{\eta}|_{\Bbb S^2}\cdot\omega
  =-\gamma\bfu_{\eta}|_{\Bbb S^2}\cdot\bfn={\gamma\over2}{\bfn}\cdot
  {\bf R}(0)[{\bfn}\cdot\Delta_{\omega}\eta]\equiv-{\gamma\over2}A_1\eta.
\eqno{(2.46)}
$$
  Define $A_0=\partial_{\bf n}(\Delta,\mbox{tr}_{\Bbb S^2})^{-1}(0,\cdot)$,
  i.e., $A_0\eta=\partial_{\bf n}\psi_\eta$ for $\eta\in h^{m+\theta}(\Bbb S^2)$,
  where $\psi_\eta$ is the solution of the following boundary value problem:
\setcounter{equation}{46}
\begin{equation}
  \Delta\psi_\eta=0 \quad\mbox{in} \;\;|x|<1,
  \quad
  \psi_\eta=\eta \quad\mbox{on}\;\;|x|=1.
\end{equation}
%(2.47)
  It is well-known that (cf. \cite{EscProk}),
\begin{equation}
  A_0\in C^\infty(h^{m+\theta}(\Bbb S^2),h^{m-1+\theta}
  (\Bbb S^2))\cap \mathcal H(h^{m+\theta}(\Bbb S^2),h^{m-1+\theta}
  (\Bbb S^2)).
\end{equation}
%(2.48)
  We rewrite
$$
  {\mathcal Q}_1(0)=-{\gamma\over 4} A_0-{\gamma\over 4}(2A_1-A_0).
$$
  In what follows we prove that
\begin{equation}
  2A_1-A_0\in L(h^{m+\theta}(\Bbb S^2),h^{m+\theta}(\Bbb S^2)).
\end{equation}
%(2.49)
  Note that if this assertion is proved, then the desired assertion
  follows.

  Since $A_0\eta=\partial_{\bf n}\psi_\eta$ and $A_1\eta=2\bfu_{\eta}|_{\Bbb S^2}
  \cdot\bfn$, we have $(2A_1-A_0)\eta=4\bfu_{\eta}|_{\Bbb S^2}\cdot\bfn-
  \partial_{\bf n}\psi_\eta$. Since $\bfu_{\eta}=-{1\over2}{\bf R}(0)
  (\Delta_{\omega}\eta\cdot{\bfn})$, by definition of the operator ${\bf R}(0)$
  we see that there exist $q_\eta\in h^{m-1+\theta}(\Bbb S^2)$ and $\zeta_\eta
  \in {\Bbb R}^6$ such that
$$
  {\mathcal L}(0)(\bfu_{\eta},q_\eta,\zeta_\eta)=\big(0,0,-{1\over2}
  \Delta_{\omega}\eta\cdot{\bfn},0,0\big).
$$
  Besides, a simple computation shows that
\begin{eqnarray*}
  {\mathcal L}(0)(\nabla\psi_\eta,0,0)&=&\Big(0,0,\bfT(\nabla\psi_\eta,0){\bfn},
  \int_{|x|<1}\nabla\psi_\eta dx,\int_{|x|<1}\nabla\psi_\eta\times x dx\Big)
  \\
  &=&\Big(0,0,2\partial_{\bfn}\nabla\psi_\eta,
  \int_{|x|=1}\eta\cdot{\bfn}dS_x,0\Big),
\end{eqnarray*}
  Since $\displaystyle 0=\Delta\psi_\eta\big|_{r=1}=
  \Big({1\over r^2}\frac{\partial^2\psi_\eta}{\partial r^2}+{2\over r}
  \frac{\partial\psi_\eta}{\partial r}+{1\over r^2}\Delta_\omega\psi_\eta
  \Big)\Big|_{r=1}=\frac{\partial^2\psi_\eta}{\partial r^2}\Big|_{r=1}+
  2\partial_{\bfn}\psi_\eta+\Delta_\omega\eta$, we have
$$
\begin{array}{rl}
  \displaystyle
  \partial_{\bfn}\nabla\psi_\eta=& \displaystyle
  {\partial\over\partial r}\nabla\psi_\eta\Big|_{r=1}=
  \frac{\partial^2\psi_\eta}{\partial r^2}\Big|_{r=1}\bfn
  +\nabla_\omega\Big(\frac{\partial\psi_\eta}{\partial r}\Big|_{r=1}\Big)
  \\ [0.3cm]
  =& \displaystyle-\Delta_\omega\eta\cdot\bfn-
  2\partial_{\bfn}\psi_\eta\cdot\bfn+\nabla_\omega(\partial_{\bfn}\psi_\eta).
\end{array}
$$
  Hence
$$
  {\mathcal L}(0)(4\bfu_\eta-\nabla\psi_\eta,4q_\eta,4\zeta_\eta)
  =\Big(0,0,4\partial_{\bfn}\psi_\eta\cdot\bfn-
  2\nabla_\omega(\partial_{\bfn}\psi_\eta),
  \int_{|x|=1}\eta\cdot{\bfn}dS_x,0\Big).
$$
  It follows that
$$
\begin{array}{rl}
  \displaystyle
  (2A_1-A_0)\eta=&\displaystyle 4\bfu_{\eta}|_{\Bbb S^2}\cdot\bfn-
  \partial_{\bf n}\psi_\eta=\displaystyle (4\bfu_{\eta}
  -\nabla\psi_\eta)|_{\Bbb S^2}\cdot\bfn \\ [0.3cm]
  =&\displaystyle
  \mbox{tr}_{\Bbb S^2}\Big\{{\mathcal J}{\mathcal L}(0)^{-1}
  \Big(0,0,4\partial_{\bfn}\psi_\eta\cdot\bfn-
  2\nabla_\omega(\partial_{\bfn}\psi_\eta),
  \int_{|x|=1}\eta\cdot{\bfn}dS_x,0\Big)\Big\}\cdot\bfn \\ [0.3cm]
  =&\displaystyle
  -2{\bf R}(0)(\nabla_\omega(\partial_{\bfn}\psi_\eta))\big|_{\Bbb S^2}\cdot\bfn
  +\mbox{tr}_{\Bbb S^2}\Big\{{\mathcal J}{\mathcal L}(0)^{-1}
  \Big(0,0,4\partial_{\bfn}\psi_\eta\cdot\bfn,
  \int_{|x|=1}\eta\cdot{\bfn}dS_x,0\Big)\Big\}\cdot\bfn \\ [0.3cm]
  \equiv& B_1\eta+B_0\eta.
\end{array}
$$
  It can be easily seen that $B_0\in L(h^{m+\theta}(\Bbb S^2),h^{m+\theta}
  (\Bbb S^2))$. Furthermore, minor changes to the proof of Lemma A.2 in
  \cite{EscProk} show that also $B_1\in L(h^{m+\theta}(\Bbb S^2),h^{m+\theta}
  (\Bbb S^2))$ (see Lemma B.1 and Corollary B.2 in Appendix B for details).
  Hence (2.49) follows. This completes the proof of Lemma 2.6.  $\qquad$$\Box$
\medskip

  Since $\mathcal H(h^{m+\theta}(\Bbb S^2),h^{m-1+\theta}(\Bbb S^2))$
  is open in $L(h^{m+\theta}(\Bbb S^2),h^{m-1+\theta}(\Bbb S^2))$,
  from Lemma 2.6 we immediately get
\medskip

  {\bf Corollary 2.7}\ \ {\em For sufficiently small $\delta$ we have
\begin{equation}
  -D{\mathcal Q}(\rho)\in\mathcal H(h^{m+\theta}(\Bbb S^2),
  h^{m-1+\theta}(\Bbb S^2)) \quad \mbox{for}\;\;\rho\in
  O_\delta^{m+\theta}(\Bbb S^2).
\end{equation}
%(2.42)
}

  By this corollary we see that, at least in a small neighborhood of the
  origin, the differential equation (2.39) is of the parabolic type in the
  sense of Amann \cite{Amann} and Lunardi \cite{Lunar}, so that the
  geometric theory for parabolic differential equations in Banach spaces
  presented in these literatures can be applied to (2.39). In the following
  sections we shall use this theory to prove Theorem 1.1.

\section{Linearization}
\setcounter{equation}{0}

  In this section we compute the spectrum of the operator $D\mathcal Q(0)$.
  Note that since $h^{m+\theta}(\Bbb S^2)$ is compactly embedded into
  $h^{m-1+\theta}(\Bbb S^2)$, by Lemma 2.6 we see that the spectrum of the
  operator $D\mathcal Q(0)$ consists of all eigenvalues.

  To compute the eigenvalues of $D\mathcal Q(0)$ we first derive
  a useful expression of this operator. Consider a perturbation of the
  radially symmetric stationary solution $(\sigma_s,\bfv_s,p_s,0)$
  (see (1.12)):
$$
  \sigma(x,t)=\sigma_s(r)+\varepsilon\phi(r,\omega,t),\quad
  \bfv(x,t)=\bfv_s(x)+\varepsilon\,\vec{\upsilon}(r,\omega,t),\quad
  p(x,t)=p_s(r)+\varepsilon\,\psi(r,\omega,t),
$$
$$
  \Omega(t)=\{x\in {\Bbb R}^3: r<1+\varepsilon\eta(\omega,t)\}
   \quad (r=|x|,\; \omega\in \Bbb S^2),
$$
  where $\varepsilon$ is a small parameter, and $\phi$, $\vec{\upsilon}
  \big(=(\upsilon_1,\upsilon_2,\upsilon_3)\big)$, $\psi$ and $\eta$ are
  new unknown functions. From \cite{FriedHu1} and \cite{FriedHu2} we see
  that the linearizations of equations (1.1)--(1.8) are respectively as
  follows:
\begin{equation}
  \Delta\phi=f'(\sigma_s)\phi
  \qquad\mbox{in}\;\;\Bbb B^3,
\end{equation}
%(3.1)
\begin{equation}
  \nabla\cdot\vec{\upsilon}=g'(\sigma_s)\phi
  \qquad\mbox{in}\;\;\Bbb B^3,
\end{equation}
%(3.2)
\begin{equation}
  -\Delta\vec{\upsilon}+\nabla\psi-
  {1\over3}\nabla(\nabla\cdot\vec{\upsilon}\,)=0
  \qquad\mbox{in}\;\;\Bbb B^3,
\end{equation}
%(3.3)
\begin{equation}
  \phi=-\sigma'_s(1)\eta
  \qquad\mbox{on}\;\;\Bbb S^2,
\end{equation}
%(3.4)
\begin{equation}
  \bfT(\vec{\upsilon},\psi){\bf n}=-2g(1)
  \nabla_\omega\eta+\gamma(\eta+{1\over2}
  \Delta_\omega\eta){\bf n}
  +4g(1)\eta\bfn\qquad\mbox{on}\;\;\Bbb S^2,
\end{equation}
%(3.5)
\begin{equation}
  \partial_t\eta=\vec{\upsilon}\,\big|_{\Bbb S^2}
  \cdot\bfn+g(1)\eta\qquad\mbox{on}\;\;\Bbb S^2,
%  =\vec{\upsilon}\,\big|_{\Bbb S^2}\cdot\bfn+
%  \eta{\partial_\bfn}{v_s}\big|_{\Bbb S^2}\cdot\bfn
\end{equation}
%(3.6)
\begin{equation}
  \int_{|x|<1}\vec{\upsilon}\,dx=0,
\end{equation}
%(3.7)
\begin{equation}
  \int_{|x|<1}\vec{\upsilon}\times x\,dx=0,
\end{equation}
%(3.8)

  Similarly as before, the system (3.1)--(3.8) can be reduced into a scalar
  equation in the unknown function $\eta$ only. Indeed, given $\eta\in
  C([0,\infty),h^{m+\theta}(\Bbb S^2))$, we first solve the second-order
  elliptic equation (3.1) subject to the boundary condition (3.4) to get
  $\phi(\cdot,t)\in h^{m+\theta}(\overline{\Bbb B}^3)$ as a functional of
  $\eta$, and next substitute this solution $\phi$ into (3.2). It can be easily
  checked that (2.27) and (2.28) are satisfied by equations (3.2), (3.3),
  (3.5), (3.7) and (3.8). Thus by using Lemma 2.3 we get a unique solution
  $(\vec{\upsilon}(\cdot,t),\psi(\cdot,t))\in h^{m-1+\theta}
  (\overline{\Bbb B}^3)\times h^{m-2+\theta}(\overline{\Bbb B}^3)$ as a
  functional of $\eta$. Substituting $\vec{\upsilon}=\vec{\upsilon}(r,\omega,t)$
  obtained in this way into (3.6) and denoting
\begin{equation}
  \mathcal B_\gamma\eta=\vec{\upsilon}\,\big|_{\Bbb S^2}\cdot\bfn+g(1)\eta,
\end{equation}
%(3.9)
  we see that the system of equations (3.1)--(3.8) reduces into the scalar
  equation
\begin{equation}
  \partial_t\eta=\mathcal B_\gamma\eta.
\end{equation}
%(3.10)
  Now, since the problem (1.1)--(1.8) is equivalent to the equation (2.39) with
  $\mathcal Q(\rho)$ given by (2.38), its linearization should correspondingly
  be equivalent to the linearization of (2.39) which reads as follows:
\begin{equation}
  \partial_t\eta=D\mathcal Q(0)\eta.
\end{equation}
%(3.11)
%  Note that we have assumed that $R_s=1$ for simpicity in Section 2, we
%  have for any $\eta\in h^{m+\theta}(\Bbb S^2)$, $\eta\in \mbox{Dom}
%  (D\mathcal Q(0))$, and we easily see that (3.11) is well defined.
  Comparing (3.10) with (3.11), we get the following result:
\medskip

  {\bf Lemma 3.1}\hs {\em $D\mathcal Q(0)=\mathcal B_\gamma$.}
  $\qquad$ $\Box$
\medskip

  In the sequel we deduce the expression of $\mathcal B_\gamma$ in terms of
  Fourier expansions of functions on the sphere ${\Bbb S}^2$.

  For each $l\in \Bbb N\cup\{0\}$, let $Y_{lm}(\omega)$ $(m=-l,-l+1,\cdots,l-1,
  l)$ be a normalized orthogonal basis of the space of all spherical harmonics
  of degree $l$. Then $\{Y_{lm}(\omega):l=0,1,2,\cdots;\,m=-l,-l+1,\cdots,l-1,
  l\}$ is a normalized orthogonal basis of the scalar $L^2$-space on
  $\Bbb S^2$. As in \cite{Hill}, let $\vec{V}_{lm}(\omega)$,
  $\vec{X}_{lm}(\omega)$ and $\vec{W}_{lm}(\omega)$, where $l=0,1,2,\cdots$
  and $m=-l,-l+1,\cdots,l-1,l$, be the corresponding vector spherical harmonics.
  From \cite{Hill} we know that all these vector spherical harmonics form a
  normalized orthogonal basis of the vector $L^2$-space on $\Bbb S^2$ (see also
  Appendix A of \cite{FriedHu1} and Section 2 of \cite{FriedReit3} for this
  assertion). Besides, for every $l\in\Bbb N\cup\{0\}$ we denote
$$
  L_l={d^2\over d r^2}+{2\over r}
  {d\over d r}-{l^2+l\over r^2}.
$$
  For simplicity of the notation, we shall not write out the whole expansions of
  $\phi$, $\vec{\upsilon}$, $\psi$ and $\eta$, but instead merely consider each
  monomials in the expansions of these functions. This is reasonable because of
  the special forms of the operators appearing in the (3.1)--(3.8). Thus we put
\begin{equation*}
  \eta=Y_{lm}(\omega).
\end{equation*}
  Then it can be easily verified that the corresponding solution of (3.1) and
  (3.4) is as follows:
\begin{equation}
  \phi(r,\omega)=F_{l}(r)Y_{lm}(\omega),
\end{equation}
%(3.12)
  where $F_{l}(r)$ is the unique solution of the following problem:
\begin{equation}
  \displaystyle L_l F_{l}(r)=f'(\sigma_s(r))F_{l}(r)
  \quad\mbox{for}\;\; 0<r<1,\quad F_{l}'(0)=0,\quad
  F_{l}(1)=-\sigma'_s(1).
\end{equation}
%(3.13)
  Observe that by (3.2) and (3.3) we have
\begin{equation}
  \Delta(\psi-{4\over3}g'(\sigma_s)\phi)=0
  \qquad\mbox{in}\;\;\Bbb B^3,
\end{equation}
%(3.14)
\begin{equation}
  \Delta\vec{\upsilon}-\nabla(g'(\sigma_s)\phi)
  =\nabla(\psi-{4\over3}g'(\sigma_s)\phi)
  \qquad\mbox{in}\;\;\Bbb B^3.
\end{equation}
%(3.15)
  Thus the solution of (3.2), (3.3), (3.5), (3.7) and (3.8) has the following
  expressions:
\begin{eqnarray}
  &\psi(r,\omega)&={4\over3}g'(\sigma_s(r))\phi(r,\omega)
  +P_{lm}(r)Y_{lm}(\omega)
  \nonumber\\
  &&={4\over3}g'(\sigma_s(r))F_{l}(r)Y_{lm}(\omega)
  +P_{lm}(r)Y_{lm}(\omega),
\end{eqnarray}
%(3.16)
\begin{equation}
  \vec{\upsilon}=\bfa+\bfb\times x+v_{lm}\vec{V}_{lm}
  +x_{lm}\vec{X}_{lm}+w_{lm}\vec{W}_{lm},
\end{equation}
%(3.17)
  where $\bfa$, $\bfb$ are unknown constant vector, and $P_{lm}(r)$, $v_{lm}(r)$,
  $x_{lm}(r)$ and $w_{lm}(r)$ are unknown functions defined on $[0,1]$ such that
  $v_{lm}(0)=x_{lm}(0)=w_{lm}(0)=0$. Using some well-known formulas for vector
  spherical harmonics (see \cite{Hill} or \cite{FriedHu1}, \cite{FriedReit3}),
  we have
\begin{eqnarray}
  \nabla\phi
%  &=&\nabla(F_{l}(r)Y_{lm}(\omega))
%  \nonumber\\
%  &=&F_{l}'(r)Y_{lm}(\omega)\bfn+F_{l}(r){l\over r}
%  \sqrt{l+1\over 2l+1}\vec{V}_{lm}(\omega)+F_{l}(r){l+1\over r}
%  \sqrt{l\over 2l+1}\vec{W}_{lm}(\omega)
%  \nonumber\\
  &=&\Big[-F_l'(r)+{l\over r}F_l(r)\Big]\sqrt{l+1\over 2l+1}\vec{V}_{lm}(\omega)
  +\Big[F'_l(r)+{l+1\over r}F_l(r)\Big]\sqrt{l\over 2l+1}\vec{W}_{lm}(\omega),
\end{eqnarray}
%(3.19)
\begin{eqnarray}
  \nabla(g'(\sigma_s)\phi)
%  &=&\,\nabla(g'(\sigma_s(r))F_l(r)Y_{lm}(\omega)) \nonumber\\
  &=&\,\Big[-{\partial\over\partial r}[g'(\sigma_s(r))F_l(r)]
  +{l\over r}g'(\sigma_s(r))F_l(r)\Big]\sqrt{l+1\over 2l+1}
  \vec{V}_{lm}(\omega)\nonumber\\
  &&+\Big[{\partial\over\partial r}[g'(\sigma_s(r))F_l(r)]
  +{l+1\over r}g'(\sigma_s(r))F_l(r)\Big]\sqrt{l\over 2l+1}
  \vec{W}_{lm}(\omega)\nonumber\\
  &\equiv&\;\digamma_l^1(r)\vec{V}_{lm}(\omega)+\digamma_l^2(r)
  \vec{W}_{lm}(\omega),
\end{eqnarray}
%(3.20)
\begin{equation}
  \nabla(\psi-{4\nu\over3}g'(\sigma_s)\phi)=\sqrt{l+1\over 2l+1}
  \big[-{\partial\over\partial r}P_{lm}+{l\over r}P_{lm}\big]\vec{V}_{lm}
  +\sqrt{l\over2l+1}\big[{\partial\over\partial r}P_{lm}+
  {l+1\over r}P_{lm}\big]\vec{W}_{lm},
\end{equation}
%(3.18)
\begin{eqnarray}
  \nabla\cdot\vec{\upsilon}
%  \;&=&\;\nabla\cdot\big[v_{lm}\vec{V}_{lm}
%  +x_{lm}\vec{X}_{lm}+w_{lm}\vec{W}_{lm}\big]\nonumber\\
  &=&\,-\big[v_{lm}'(r)+{l+2\over r}v_{lm}(r)\big]\sqrt{l+1\over2l+1}
  Y_{lm}(\omega)+
  \big[w_{lm}'(r)-{l-1\over r}w_{lm}(r)\big]\sqrt{l\over2l+1}
  Y_{lm}(\omega),\nonumber\\
  &&
\end{eqnarray}
%(3.21)
  and
\begin{equation}
  \Delta\vec{\upsilon}=L_{l+1}(v_{lm}(r))\vec{V}_{lm}(\omega)
  +L_l(x_{lm}(r))\vec{X}_{lm}(\omega)+L_{l-1}(w_{lm}(r))
  \vec{W}_{lm}(\omega).
\end{equation}
%(3.22)
  By (3.15), (3.18), (3.20) and (3.22), we have
\begin{equation}
  \sqrt{l+1\over 2l+1}\Big[-P_{lm}'(r)+{l\over r}P_{lm}(r)\Big]
  =L_{l+1}(v_{lm}(r))-\digamma_l^1(r),
\end{equation}
%(3.23)
\begin{equation}
  L_l(x_{lm}(r))=0.
\end{equation}
%(3.24)
\begin{equation}
  \sqrt{l\over 2l+1}\Big[P_{lm}'(r)+{l+1\over r}P_{lm}(r)\Big]
  =L_{l-1}(w_{lm}(r))-\digamma_l^2(r).
\end{equation}
%(3.25)
  By (3.14) we have
\begin{equation}
  L_l(P_{lm}(r))=0.
\end{equation}
%(3.26)
  By (3.2) and (3.21) we have
\begin{equation}
  \sqrt{l\over2l+1}\Big[w_{lm}'(r)-{l-1\over r}w_{lm}(r)\Big]-
  \sqrt{l+1\over2l+1}\Big[v_{lm}'(r)+{l+2\over r}v_{lm}(r)\Big]
  =g'(\sigma_s(r))F_l(r).
\end{equation}
%(3.27)
  Solving the ODE problem (3.23)--(3.27), we get
\begin{equation}
  P_{lm}(r)=2(2l+3)A_1r^l, \qquad
  x_{lm}(r)=B_1r^l,
\end{equation}
%(3.28)
\begin{eqnarray}
  v_{lm}(r)\;&=&\;\sqrt{l+1\over2l+1}{2l\over l+1}A_1r^{l+1}-
  {r\over(2l+3)}\sqrt{l+1\over2l+1}g'(\sigma_s(r))F_l(r)
  -{r^{-l-2}\over2l+3}\int_0^r s^{l+3}\digamma_l^1(s)\,ds
  \nonumber
  \\
  &\equiv&\;\sqrt{l+1\over2l+1}{2l\over l+1}A_1r^{l+1}
  -\tilde v_l(r),
\end{eqnarray}
%(3.29)
  and
\begin{eqnarray}
  w_{lm}(r)\;&=&\;C_1r^{l-1}+\sqrt{l\over2l+1}(2l+3)A_1
  r^{l+1}\nonumber\\
  &&-{r\over(2l-1)}\sqrt{l\over2l+1}g'(\sigma_s(r))
  F_l(r)-{r^{l-1}\over2l-1}\int_r^{R_s}s^{-l+2}
  \digamma_l^2(s)\,ds\nonumber\\
  &\equiv&\;C_1r^{l-1}+\sqrt{l\over2l+1}(2l+3)A_1r^{l+1}
  -\tilde w_l(r),
\end{eqnarray}
%(3.30)
  where $A_1$, $B_1$ and $C_1$ are constants.

  Next we consider the boundary condition (3.5). Again by using
  some well-known properties of vector spherical harmonics
  (see \cite{Hill} or \cite{FriedHu1}, \cite{FriedReit3}),
  we can rewrite $\vec{\upsilon}$ in (3.17) as follows
\begin{equation}
  \vec{\upsilon}(r,\omega)=\bfa+\bfb\times x+H_{l1}(r)Y_{lm}(\omega)
  \,\omega+H_{l2}(r)\nabla_\omega Y_{lm}(\omega),
\end{equation}
%(3.31)
  where
$$
  H_{l1}(r)=-\sqrt{l+1\over2l+1}v_{lm}(r)+\sqrt{l\over2l+1}w_{lm}(r),
  \qquad
  H_{l2}(r)={v_{lm}(r)\over\sqrt{(l+1)(2l+1)}}+{w_{lm}(r)\over\sqrt{l(2l+1)}}.
$$
  Thus
\begin{eqnarray}
  \bfT(\vec{\upsilon},\psi){\bf n}
  \;&=&\;{2\over3}g'(1)\sigma'_s(1)Y_{lm}(\omega)\,\omega
  +\big[2H_{l1}'(1)Y_{lm}(\omega)-\psi(1,\omega)\big]\omega
  \nonumber\\
  &&+\big[{H_{l1}(1)-H_{l2}(1)}+H_{l2}'(1)\big]\nabla_\omega
  Y_{lm}(\omega).
\end{eqnarray}
%(3.32)
  Note that by (3.16) and (3.28) we have
$$
  \psi(1,\omega)=Y_{lm}(\omega)\big[-{4\over3}g'(1)
  \sigma_s'(1)+2(2l+3) A_1\big].
$$
  Substituting this expression into (3.32) we get
\begin{eqnarray}
  \bfT(\vec{\upsilon},\psi){\bf n}
  \;&=&\;\big[{2\over3}g'(1)\sigma'_s(1)+2H_{l1}'(1)+{4\over3}g'(1)
  \sigma_s'(1)-2(2l+3) A_1\big]Y_{lm}(\omega)\,\omega
  \nonumber\\
  &&+\big[{H_{l1}(1)-H_{l2}(1)}+H_{l2}'(1)\big]\nabla_\omega
  Y_{lm}(\omega).
\end{eqnarray}
%(3.33)
  On the other hand, putting $\eta=Y_{lm}(\omega)$ in (3.5) and using the
  well-known relation $\Delta_\omega Y_{lm}(\omega)=-l(l+1)Y_{lm}(\omega)$,
  we get
\begin{equation}
  \bfT(\vec{\upsilon},\psi){\bf n}=-2g(1)\nabla_\omega Y_{lm}(\omega)
  +\big[\gamma(1-{l(l+1)\over2})+4g(1)\big]Y_{lm}(\omega)\,\omega.
\end{equation}
%(3.34)
  Since $\nabla_\omega Y_{lm}(\omega)$ and $Y_{lm}(\omega)\,\omega$ are
  mutually orthogonal, by comparing their coefficients in (3.33) and (3.34)
  and using the relations
$$
  H_{l1}(1)=-\sqrt{l+1\over2l+1}v_{lm}(1)+\sqrt{l\over2l+1}w_{lm}(1),
  \qquad
  H_{l2}(1)={v_{lm}(1)\over\sqrt{(l+1)(2l+1)}}+{w_{lm}(1)\over\sqrt{l(2l+1)}},
$$
$$
  H_{l1}'(1)=-\sqrt{l+1\over2l+1}v_{lm}'(1)+\sqrt{l\over2l+1}w_{lm}'(1),
  \qquad
  H_{l2}'(1)={v_{lm}'(1)\over\sqrt{(l+1)(2l+1)}}+{w_{lm}'(1)\over\sqrt{l(2l+1)}},
$$
  we obtain
\begin{eqnarray}
  -\sqrt{l+1\over2l+1} v_{lm}'(1)+\sqrt{l\over2l+1}w_{lm}'(1)
  =\;{\gamma\over4}(2-l^2-l)+2g(1)-g'(1)
  \sigma_s'(1)+(2l+3)A_1,
\nonumber\\
\end{eqnarray}
%(3.35)
  and
\begin{equation}
  {1\over\sqrt{2l+1}}\Big[-{l+2\over\sqrt{l+1}}v_{lm}(1)+
  {l-1\over \sqrt{l}}w_{lm}(1)\Big]+{v_{lm}'(1)\over\sqrt{(l+1)(2l+1)}}
  +{w_{lm}'(1)\over\sqrt{l(2l+1)}}=-2g(1).
\end{equation}
%(3.36)

  We now proceed to consider the equation (3.6). By (3.9) and (3.29)--(3.31)
  we have
\begin{eqnarray}
  \mathcal B_\gamma Y_{lm}(\omega)\;&=&\;\vec{\upsilon}\,
  \big|_{\Bbb S^2}\cdot{\bf n}+g(1)Y_{lm}(\omega)=\bfa\cdot\omega
  +[H_{l1}(1)+g(1)]Y_{lm}(\omega)\nonumber\\
  &=&\bfa\cdot\omega+Y_{lm}(\omega)\Big[g(1)-\sqrt{l+1\over2l+1}v_{lm}(1)
  +\sqrt{l\over2l+1}w_{lm}(1)\Big]
  \nonumber\\
  &=&\bfa\cdot\omega+Y_{lm}(\omega)\Big[g(1)+l(A_1+\tilde C_1)
  +\sqrt{l+1\over2l+1}\tilde v_l(1)
  -\sqrt{l\over2l+1}\tilde w_l(1)\Big],
\end{eqnarray}
%(3.37)
  where $\tilde C_1=C_1/\sqrt{l(2l+1)}$. Thanks to the constrain condition
  (3.7) we see that
\begin{equation}
  \bfa=-{3\over 4\pi}\int_{|x|<1}\{H_{l1}(r)Y_{lm}(\omega)\,\omega+
  H_{l2}(r)\nabla_\omega Y_{lm}(\omega)\}\,dx=0 \quad
  \mbox{for}\;\;l\in\Bbb N,\;l\neq 1
\end{equation}
%(3.38)
  (cf. (5.8) in \cite{FriedHu1}). To compute $A_1+\tilde C_1$ we
  substitute (3.29) and (3.30) into (3.35) and (3.36), which gives
\begin{eqnarray}
  \sqrt{l\over2l+1}\tilde w_l'(1)-\sqrt{l+1\over2l+1}\tilde v_l'(1)
  &=&-{\gamma\over4}(2-l^2-l)-2g(1)+g'(1)\sigma_s'(1)  \nonumber\\
  &&+\tilde C_1(l^2-l)+(l^2-l-3)A_1,
\end{eqnarray}
%(3.39)
  and
\begin{eqnarray}
  &&\sqrt{l\over2l+1}\big[{l-1\over l}\tilde w_l(1)+{1\over l}\tilde
  w_l'(1)\big]-\sqrt{l+1\over2l+1}\big[{l+2\over l+1}\tilde v_l(1)
  -{1\over l+1}\tilde v_l'(1)\big]\nonumber
  \\
  &=&\;2g(1)+2\tilde C_1(l-1)+{2l^2+4l\over l+1}A_1.
\end{eqnarray}
%(3.40)
  $(3.39)\times2(2l+1)+(3.40)\times3(l+1)$ yields:
\begin{eqnarray}
  A_1+\tilde C_1\;&=&\;{1\over2(l-1)(2l^2+4l+3)}\Big\{
  {\gamma\over2}(1-l)(2l^2+5l+2)-(4l+2)g'(1)\sigma_s'(1)
  \nonumber\\
  &&+(2l-2)g(1)+{4l^2+5l+3\over\sqrt{l(2l+1)}}\tilde w_l'(1)
  -\sqrt{l+1\over2l+1}(4l-1)\tilde v_l'(1)
  \nonumber\\
  &&+{3(l^2-1)\over\sqrt{l(2l+1)}}\tilde w_l(1)-3(l+2)\sqrt
  {l+1\over2l+1}\tilde v_l(1)\Big\}.
\end{eqnarray}
%(3.41)
  By the definitions of $\digamma_l^i$ $(i=1,2)$, $\tilde v_l$
  and $\tilde w_l$ respectively in (3.20), (3.29) and (3.30),
  and by straightforward calculation we easily have
\begin{equation}
  \tilde v_l(1)=\sqrt{l+1\over2l+1}\int_0^1
  g'(\sigma_s(r))F_l(r)r^{l+2}\,dr,
\end{equation}
%(3.42)
\begin{equation}
  \tilde w_l(1)=-\sqrt{l\over2l+1}{g'(1)\over2l-1}\sigma_s'(1),
\end{equation}
%(3.43)
\begin{equation}
  \tilde v_l'(1)=-\sqrt{l+1\over2l+1}g'(1)\sigma_s'(1)
  -(l+2)\sqrt{l+1\over2l+1}\int_0^1g'(\sigma_s(r))F_l(r)r^{l+2}\,dr,
\end{equation}
%(3.44)
\begin{equation}
  \tilde w_l'(1)=\sqrt{l\over2l+1}{l\over2l-1}g'(1)\sigma_s'(1).
\end{equation}
%(3.45)
  Substituting (3.38) and (3.41) into (3.37) and using (3.42)--(3.45), we see
  that, for $l\neq1$,
\begin{eqnarray}
  \mathcal B_\gamma Y_{lm}(\omega)&=&{1\over2l^2+4l+3}
  \Big\{g(1)(2l+3)(l+1)-{\gamma\over4}l(2l+1)(l+2)
  \nonumber\\
  &&+(2l+3)(l+1)\int_0^1g'(\sigma_s(r))
  F_l(r)r^{l+2}\,dr\Big\}Y_{lm}(\omega).
\end{eqnarray}
%(3.46)
  We define, for $l\geq 2$,
\begin{equation}
  \gamma_l={4(2l+3)(l+1)\over l(l+2)(2l+1)}\Big[
  g(1)+\int_0^1g'(\sigma_s(r))F_l(r)r^{l+2}\,dr\Big],
\end{equation}
%(3.47)
\begin{equation}
  \alpha_l(\gamma)=-{l(l+2)(2l+1)\over 4(2l^2+4l+3)}(\gamma-\gamma_l).
\end{equation}
%(3.48)
  Then in case $l\geq 2$ (3.46) can be rewritten as follows:
$$
  \mathcal B_\gamma Y_{lm}(\omega)=\alpha_l(\gamma)Y_{lm}(\omega).
$$
  In the case $l=0$ we have, directly from (3.46), that (3.49) also holds
  with
\begin{equation}
  \alpha_0\equiv\alpha_0(\gamma)=g(1)+\int_0^1g'(\sigma_s(r))F_0(r)r^2\,dr.
\end{equation}
%(3.49)

  Finally, we consider the case $l=1$. Since the problem (1.1)--(1.8) is
  translation invariant, by some similar argument as those in \cite{Cui2}
  and \cite{FriedHu1} we see that $\mathcal B_\gamma \eta=0$ for any sphere
  harmonics of degree $1$. In particular, we have
$$
  \mathcal B_\gamma Y_{1m}(\omega)=0, \quad m=-1,0,1.
$$

  In summary, we have proved the following result:
\medskip

  {\bf Lemma 3.2} \hs {\em $D\mathcal Q(0)=\mathcal B_\gamma $ is a Fourier
  multiplication operator having the following expression: For any
  $\eta\in C^\infty({\Bbb S}^2)$ with Fourier expansion $\eta(\omega)=
  \sum_{l=0}^\infty\sum_{m=-l}^l c_{lm}Y_{lm}(\omega)$, we have
\begin{equation}
  \mathcal B_\gamma\eta(\omega)=\alpha_0c_{00}Y_{00}+
  \sum_{l=2}^\infty\sum_{m=-l}^l\alpha_l(\gamma)c_{lm}
  Y_{lm}(\omega),
\end{equation}
%(3.50)
  where $\alpha_0$ and $\alpha_l(\gamma)$ defined in $(3.49)$ and $(3.48)$,
  respectively.} $\qquad$$\Box$
\medskip

  As usual, for a given closed linear operator $B$ in a Banach space $X$, we
  denote by $\rho(B)$ and $\sigma(B)$ the resolvent set and the spectrum of $B$,
  respectively. The set of all eigenvalues of $B$ is denoted by $\sigma_p(B)$.
  As mentioned in the beginning of this section, we have $\sigma(D\mathcal Q(0))
  =\sigma_p(D\mathcal Q(0))$. Hence, from Lemma 3.2 we immediately obtain the
  following result:
\medskip

  {\bf Lemma 3.3}\hs {\em The spectrum of $D\mathcal Q(0)=\mathcal B_\gamma$ is
  given by
$$
  \sigma(\mathcal B_\gamma)=\{\alpha_0,0\}\cup\{\alpha_l(\gamma):l=2,3,4,\cdots\}.
$$
  Moreover, the multiplicity of the eigenvalue $0$ is $3$.} $\qquad\Box$
\medskip

  The next result shows some useful properties of $\alpha_0$ and $\gamma_l$
  ($l\ge 2$):
\medskip

  {\bf Lemma 3.4}\hs {\em We have the following assertions:

  $(i)$ $\alpha_0<0$.

  $(ii)$ $\gamma_l>0$ for all $l\ge2$, and $\lim_{l\to\infty}\gamma_l=0$.

  $(iii)$ There exists an integer $l^*\ge2$ such that
  $\gamma_{l+1}<\gamma_l$ for all $l\ge l^*$. $\qquad$$\Box$}
\medskip

  {\em Proof:}\ \ By Assumption $(A1)$ we have $f'>0$. Thus, by the
  maximum principle we see that $F_0(r)\leq0$. Furthermore, since
  $u(r)=-\sigma_s'(r)$ is a solution of the problem
\begin{equation}
  u''(r)+{2\over r}u'(r)=f'(\sigma_s(r))u(r)+{2\over r^2}u
  \quad\mbox{for}\;\;0<r<1,\quad u(0)=0, \quad u(1)=-\sigma_s'(1),
\end{equation}
%(3.51)
  by comparison we easily get $F_0(r)\le -\sigma_s'(r)$. Thus, since $g'>0$
  (by Assumption $(A2)$), we have
$$
  \int_0^1g'(\sigma_s(r))F_0(r)r^2\,dr\le-\int_0^1g'(\sigma_s(r))
  \sigma_s'(r)r^2\,dr=-g(1)+2\int_0^1g(\sigma_s(r))r\,dr.
$$
  Hence,
$$
  \alpha_0=g(1)+\int_0^1g'(\sigma_s(r))F_0(r)r^2\,dr\le
  2\int_0^1g(\sigma_s(r))r\,dr<0.
$$
  The last inequality follows from the facts that $g'>0$ and $\displaystyle
  \int_0^1 g(\sigma_s(r))r^2\,dr=0$ (by $(A.3)$ and $(A.4)$ in Appendix A).
  This proves $(i)$.

  Next, from (3.13) we have, for any $l\ge2$, that
$$
  F_l''(r)+{2\over r}F_l'(r)-{2\over r^2}F_l(r)-f'(\sigma_s(r))F_l(r)
  ={l^2+l-2\over r^2}F_l(r)\le0
$$
  Since $f'>0$ and $u(r)=-\sigma_s'(r)$ satisfies (3.51), by comparison we get
  $-\sigma_s'(r)<F_l(r)<0$. Hence
\begin{eqnarray}
  g(1)+\int_0^1g'(\sigma_s(r))F_l(r)r^{l+2}\,dr
  \;&>&\; g(1)-\int_0^1g'(\sigma_s(r))\sigma_s'(r)r^{l+2}\,dr
  \nonumber\\
  &>&\;g(1)-\int_0^1g'(\sigma_s(r))\sigma_s'(r)r^3\,dr
  \nonumber\\
  &=&\;3\int_0^1g(\sigma_s(r))r^2\,dr=0,
\end{eqnarray}
%(3.52)
  so that $\gamma_l>0$ for all $l\ge2$. Moreover, since $|F_l(r)|\le
  \sigma_s'(1)$ and $0<g'(\sigma_s(r))\leq g'(1)$, we have
\begin{eqnarray*}
  \gamma_l\;&=&\;{4(2l+3)(l+1)\over l(l+2)(2l+1)}\Big[g(1)+\int_0^1
  g'(\sigma_s(r))F_l(r)r^{l+2}\,dr\Big]
%  \\
%  &\le&\;
 \leq 8l^{-1}\big[g(1)+{g'(1)\sigma_s'(1)\over l+3}\big].
\end{eqnarray*}
  Hence $\lim_{l\to\infty}\gamma_l=0$. This proves $(ii)$.

  Finally, by direct computation we have
$$
  \gamma_{l+1}-\gamma_l=-4g(1)(1+o(1))l^{-2}\quad\mbox{as}
  \;\;l\to\infty.
$$
  From this fact the assertion $(iii)$ immediately follows. The proof is
  complete. $\qquad$ $\Box$
\medskip

  By virtue of the the assertion $(ii)$ of the above lemma, we introduce
\begin{equation}
  \gamma_*=\max_{l\geq 2}\gamma_l.
\end{equation}
%(3.53)
  Note that Lemma 3.4 ensures that $0<\gamma_*<\infty$.

  It is interesting to compare the threshold number $\gamma_*$ defined
  above with the corresponding threshold number for the porous
  medium structured tumor model obtained by Cui and Escher \cite{CuiEsc2},
  which we denote by $\tilde\gamma_*$. Recall that $\tilde\gamma_*=
  \max_{l\geq 2}\tilde\gamma_l$, where, for the case $R_s=1$ and
  $\bar{\sigma}=1$,
\begin{equation}
  \tilde\gamma_l={2\over l(l-1)(l+2)}\Big[g(1)+\int_0^1g'(1)F_l(r)r^{l+2}
  \,dr\Big]\qquad\mbox{for}\;\;l\ge2,\;\l\in\Bbb N.
\end{equation}
%(3.54)
  From Lemma 3.2 of \cite{CuiEsc2} we know that $\{\tilde\gamma_l\}_{l\ge2}$
  has the same properties as $\{\gamma_l\}_{l\ge2}$ presented in Lemma 3.4.
\medskip

  {\bf Lemma 3.5}\ \ {\em $\tilde\gamma_l<\gamma_l$ for all $l\ge2$,
  $l\in\Bbb N$, so that $\tilde\gamma_*<\gamma_*$.}
\medskip

  {\em Proof:}\ \  For any $l\ge2$ we have, by (3.47), (3.52) and (3.54), that
\begin{eqnarray*}
  \tilde\gamma_l-\gamma_l\;&=&\;\Big[{2\over l(l-1)(l+2)}-{4(2l+3)
  (l+1)\over l(l+2)(2l+1)}\Big]\Big[g(1)+\int_0^1g'(1)F_l(r)r^{l+2}
  \,dr\Big]\\
  &=&\;{2\over l(l+2)}\Big[{1\over l}-{2(2l+3)(l+1)\over 2l+1}\Big]
  \Big[g(1)+\int_0^1g'(1)F_l(r)r^{l+2}\,dr\Big]
  \\
  &<&\;{2\over l(l+2)}(1-1)\Big[g(1)+\int_0^1g'(1)F_l(r)
  r^{l+2}\,dr\Big]=0.
\end{eqnarray*}
  This completes the proof. $\qquad$ $\Box$
\medskip

\section{The proof of Theorem 1.1}
\setcounter{equation}{0}

  In this section we give the proof of Theorem 1.1. Note that since $0\in
  \sigma(D{\mathcal Q}(0))$, the standard linearized stability theorem for
  parabolic differential equations in Banach spaces cannot be applied to
  treat (2.39), and we have to employ the method of center manifold analysis.
  We shall construct a locally invariant center manifold, which consists only
  of equilibria, and show that this manifold attracts nearby transient
  solutions at an exponential rate. Similar method was applied in
  \cite{CuiEsc1}, \cite{EscSim} and \cite{WuCui}.
\medskip

  {\bf Proof of Theorem 1.2}: \hs We fulfill this proof in four steps.

  (i) By the definition of $\gamma_*$ we see that for any $\gamma
  >\gamma_*$, $\alpha_l(\gamma)<0$ for all $l\ge2$. Besides,
  by (3.50) we see that $0$ is an eigenvalue of geometric multiplicity 3,
  and the kernel of $D\mathcal Q(0)=\mathcal B_\gamma$ is the
  space $X_c:=\mbox{span}\{Y_{1m}; m=-1,0,1\}$. Let $X_c^\bot$ be the
  orthogonal complement of $X_c$ in $L^2(\Bbb S^2)$, and for fixed $m\ge3$
  we denote $h_s^{m+\theta}(\Bbb S^2)=h^{m+\theta}(\Bbb S^2)\cap X_c^\bot$.
  Then we have
$$
  h^{m+\theta}(\Bbb S^2)=h_s^{m+\theta}(\Bbb S^2)\oplus X_c.
$$
  This decomposition induces two projection operators $\pi^c$ and $\pi^s$ which
  map $h^{m+\theta}(\Bbb S^2)$ onto $X_c$ and $h_s^{m+\theta}(\Bbb S^2)$,
  respectively. From Lemma 3.2 we know that $\mathcal B_\gamma$ commutes
  with them.

  (ii) Let $M(\eta)=\mathcal Q(\eta)-\mathcal B_\gamma\eta$.
  Then the equation (2.39) can be rewritten as follows:
\begin{equation}
  \partial_t \eta=\mathcal B_\gamma\eta+M(\eta)\;\;\quad\mbox{for}\;\;
  t>0,\quad\mbox{and}\quad \eta(0)=\eta_0.
\end{equation}
%(4.1)
  The little H\"{o}lder spaces have the following well-known
  interpolation property
$$
  (h^{\sigma_0}(\Bbb S^2),h^{\sigma_1}(\Bbb S^2))_\vartheta=
  h^{(1-\vartheta)\sigma_0+\vartheta\sigma_1}(\Bbb S^2), \qquad
  \mbox{if}\;\;(1-\vartheta)\sigma_0+\vartheta\sigma_1\notin \Bbb Z
$$
  where $0<\vartheta<1$ and $(\cdot,\cdot)_\vartheta$ denotes the
  continuous interpolation of Da Prato and Grisvard (see \cite{Lunar}).
  By Lemma 3.1 we know that $\mathcal B_\gamma=D\mathcal Q(0)$ generates
  a strongly continuous analytic semigroup on $h^{m-1+\theta}({\Bbb S}^2)$
  with domain $h^{m+\theta}({\Bbb S}^2)$. Thus by Propositions 6.2, 6.4
  and Theorem 6.5 in \cite{EscSim} we conclude that there exists an open
  neighborhood ${\mathcal O}$ of the origin in $X_c$ and a mapping
$$
  \mbox{\tensy C}\in C^m({\mathcal O},h^{m+\theta}_s(\Bbb S^2))
  \quad \mbox{with}\;\; \mbox{\tensy C}(0)=0, \quad \partial
  \mbox{\tensy C}(0)=0,
$$
  such that the 3-dimensional submanifold {\tensy M}$_c:=\,\mbox{graph}\,
  (\mbox{\tensy C})$ of $h^{m+\theta}(\Bbb S^2)$ is a locally invariant
  and stable manifold for the evolution equation (4.1). Note that
  $\mathcal M_c$ consists only of radial equilibria, i.e. $\mathcal M_c$
  is the set of all spheres of radius 1 with centers sufficiently
  close to 0. Furthermore, by the above-mentioned results of \cite{EscSim}
  we know that $\mathcal M_c$ attracts at an exponential rate all small
  global solutions of (4.1) in $h^{m+\theta}(\Bbb S^2)$. More precisely,
  there exists $\epsln>0$ such that the solution to (4.1) exists globally
  for any $\eta_0$ with $\|\eta_0\|_{h^{m+\theta}(\Bbb S^2)}\le\epsln$,
  and, moreover, there exist $c>0$, $K>0$ and a unique $z_0=z_0(\eta_0)
  \in{\mathcal O}$ such that for any $t\ge0$ there holds
\begin{equation}
  \|(\pi^c\eta(t),\pi^s\eta(t))-(z_0, \mbox{\tensy C}(z_0))\|
  _{h^{m+\theta}(\Bbb S^2)}\le K\exp(-ct)\|\pi^s\eta_0-\mbox
  {\tensy C}(\pi^c\eta_0)\|_{h^{m+\theta}(\Bbb S^2)}.
\end{equation}
%(4.2)

  (iii) Now let $\eta_0\in h^{m+\theta}(\Bbb S^2)$ be given and
  $\|\eta_0\|_{C^{m+\theta}(\Bbb S^2)}\le\epsln$. Then the
  solution of the equation (4.1)
$
  \eta\in C([0,\infty),h^{m+\theta}(\Bbb S^2))\cap C^1
  ((0,\infty),h^{m+\theta-1}(\Bbb S^2)),
$
  and it satisfies (4.1). By Lemma 2.1 and Lemma 2.5, it follows
  that the problem (1.1)--(1.9) has a global-in-time solution
  $(\sigma(\cdot,t),v(\cdot,t),p(\cdot,t),\Omega(t))$, where
  $\Omega(t)=\{x\in \Bbb R^3:\; x=r\omega,\;0\le r<1+\eta(\omega,t),
  \; \omega\in\Bbb S^2\}$. Since $\mathcal M_c$ is the set of
  equilibrium solutions which are sufficiently close to $\Bbb S^2$,
  there exists a $x_0\in \Bbb R^3$ such that $(z_0,\mbox{\tensy C}(z_0))
  =\eta_{[x_0]}$, where $\eta_{[x_0]}$ is the distance function on
  $\Bbb S^2$ introduced in Section 1. Then (4.2) implies that
\begin{equation}
  \|\eta(\cdot,t)-\eta_{[x_0]}\|_{h^{m+\theta}(\Bbb S^2)}\le
  K\exp(-ct) \qquad \mbox{for any}\;\; t\ge0.
\end{equation}
%(4.3)
  By Lemmas 2.2--2.4 we have
\begin{equation}
   \sigma(\cdot,t)=\mathcal R(\eta(t))\circ\Phi_{\eta(t)}^{-1},
   \qquad \bfv(\cdot,t)=\widetilde\bfv(\eta(t))\circ\Phi_{\eta(t)}^{-1},
   \qquad p(\cdot,t)=\widetilde p(\eta(t))\circ\Phi_{\eta(t)}^{-1},
\end{equation}
%(4.4)
  Recalling the definition of $(\sigma_{[x_0]},\bfv_{[x_0]},p_{[x_0]},
  \Omega_{[x_0]})$, we have
\begin{equation}
   \sigma_{[x_0]}=\mathcal R(\eta_{[x_0]})\circ\Phi_{\eta_{[x_0]}}^{-1},
   \qquad v_{[x_0]}=\widetilde \bfv(\eta_{[x_0]})\circ\Phi_{\eta_{[x_0]}}^{-1},
   \qquad p_{[x_0]}=\widetilde p(\eta_{[x_0]})\circ\Phi_{\eta_{[x_0]}}^{-1},
\end{equation}
%(4.5)
  The explicit construction of $\Phi_\eta$, (2.21) and the mean value theorem
  immediately imply that there is a positive constant $C$ such that
  for any $t\ge0$,
$$
  \|\mathcal R(\eta(t))-\mathcal R(\eta_{[x_0]})\|
  _{C^{m+\theta}(\overline{\Bbb B}^3)}\le C \|\eta(t)-\eta_{[x_0]}\|
  _{C^{m+\theta}(\overline{\Bbb B}^3)}.
$$
  Then using (4.4), (4.5) we have
\begin{equation}
  \|\sigma(\cdot,t)-\sigma_{[x_0]})\|_{C^{m+\theta}(\bar\Omega(t))}
  \le C\|\eta(t)-\eta_{[x_0]}\|_{C^{m+\theta}(\bar\Omega(t))},
\end{equation}
%(4.6)
  for any $t\ge0$. Similarly, by Lemma 2.3 and Lemma 2.4 we also have
\begin{equation}
  \|\bfv(\cdot,t)-\bfv_{[x_0]})\|_{C^{m-1+\theta}(\bar\Omega(t))}
  \le C\|\eta(t)-\eta_{[x_0]}\|_{C^{m-1+\theta}(\bar\Omega(t))},
\end{equation}
%(5.7)
\begin{equation}
  \|p(\cdot,t)-p_{[x_0]})\|_{C^{m-2+\theta}(\bar\Omega(t))}
  \le C\|\eta(t)-\eta_{[x_0]}\|_{C^{m-2+\theta}(\bar\Omega(t))},
\end{equation}
%(5.8)
  for any $t\ge0$. Combining (4.2), (4.3) and (4.6)--(4.8), we see that
  (1.14) holds.

  (iv) Finally, if $0<\gamma<\gamma_*$ then by Lemma 3.3 we see that
  $\sigma(\mathcal B_\gamma)\cap\{\lambda\in\Bbb C: {\rm Re}\lambda>0\}$
  is not empty. It follows from Theorem 9.1.3 in \cite{Lunar} that the
  zero equilibrium of (4.1) is unstable. This completes the proof of
  Theorem 1.1. $\qquad$ $\square$
\medskip

\section*{Appendix A:\ \ Radially symmetric stationary solution}

\hspace{2em}
  {\bf Theorem A}\ \ {\em Under Assumptions $(A1)$--$(A3)$, the problem
  $(1.1)$--$(1.8)$ has a unique radially symmetric stationary solution
  $(\sigma_s, \bfv_s, p_s, \Omega_s)$ with components having expressions
  in $(1.12)$.}
\medskip

  {\em Proof:} \ \
  It suffices to consider the following problem:
$$
   \sigma_s''(r)+{2\over r}\sigma_s'(r)=f(\sigma_s(r)) \quad
   \mbox{for} \hs 0<r\le R_s,
\eqno{(A.1)}
$$
$$
  \sigma_s'(0)=0,
  \qquad
  \sigma_s(R_s)=1,
\eqno{(A.2)}
$$
$$
  v_s'(r)+{2\over r}v_s(r)=g(\sigma_s(r)) \quad \mbox{for} \hs 0<r\le R_s,
\eqno{(A.3)}
$$
$$
  v_s(0)=0,
  \qquad
  v_s(R_s)=0.
\eqno{(A.4)}
$$
$$
  \big[{d^2\over d r^2}+{2\over r}{d\over d r}\big]
  \big[p_s(r)-{4\over3}g(\sigma_s(r))\big]=0 \quad
  \mbox{for} \hs 0<r\le R_s,
\eqno{(A.5)}
$$
$$
  p'_s(0)=0, \qquad p_s(R_s)={\gamma\over
  R_s}+{4\over3}g(1).
\eqno{(A.6)}
$$
  Indeed, it is straightforward to verify that for a radially symmetric
  stationary solution the equations (1.1) and (1.2) have respectively
  the forms $(A.1)$ and $(A.3)$, the boundary condition (1.4) has the form
  of the second equation in $(A.2)$, and (1.6) becomes the second equation in
  $(A.4)$. Next, taking divergence in both sides of (1.3) and using (1.2), we
  see that $(A.5)$ holds, by which we have
\begin{eqnarray*}
  \bfT(\bfv_s,p_s){\bfn}\big|_{r=R_s}\;&=&\;\big[2 v_s'(R_s)-p_s(R_s)-
  {2\over3}g(1)\big]{\bf n}
  \\
  &=&\;[2 g(1)-p_s(R_s)-{2\over3}g(1)]{\bf n}
  \\
  &=&\;[{4\over3} g(1)-p_s(R_s)]{\bf n},
  \qquad(\mbox{by}\;\;(A.3),\;(A.4))
\end{eqnarray*}
  where ${\bf n}(x)=x/|x|$ for $x\in\Bbb R^3\backslash\{0\}$. Hence, the second
  equation in $(A.6)$ follows from the boundary condition (1.5). Finally, the
  first equations in $(A.2)$, $(A.4)$ and $(A.6)$ are imposed to rule out
  possible solutions possessing singularities at $r=0$ for the problem without
  these equations, which are not meaningful solutions of (1.1)--(1.8).

  From \cite{Cui0} we know that under Assumptions $(A1)$--$(A3)$, the problem
  $(A.1)$--$(A.4)$ has a unique solution $(\sigma_s(r), v_s(r),R_s)$. Besides,
  using $(A.5)$ we immediately see that the function
$$
  p_s(r)={\gamma\over R_s}+{4\over3}g(\sigma_s(r))
$$
  solves $(A.5)$ and $(A.6)$. Since the solution of $(A.1)$--$(A.6)$ is
  obviously unique, we see that the desired assertion follows. $\qquad$$\Box$
\medskip

\section*{Appendix B:\ \ Boundedness of the operator $B_1$}

\hspace{2em}
  {\bf Lemma B.1}\ \ {\em Let $m\in\Bbb N$, $m\ge3$, and $0<\theta<1$.
  Let $\bfh\in(h^{m-2+\theta}(\Bbb S^2))^3$ satisfy $\bfh\cdot{\bf n}=0$
  and let $\bfv={\mathbf R}(0)\bfh$. Then $\mbox{\rm tr}_{\Bbb S^2}(\bfv)
  \cdot {\bf n}\in h^{m+\theta}(\Bbb S^2)$, and there exists a constant
  $C$ independent of $\bfh$ such that
$$
  \|\mbox{\rm tr}_{\Bbb S^2}(\bfv)\cdot{\bf n}\|_{h^{m+\theta}(\Bbb S^2)}\le
  C\|\bfh\|_{(h^{m-2+\theta}(\Bbb S^2))^3}.
$$
}

  {\em Proof:}\ \ Fix a function $d\in BUC^\infty({\Bbb R}^3)$ such
  that
$$
  \mbox{tr}_{\Bbb S^2}(d)=0,\qquad
  \mbox{tr}_{\Bbb S^2}(\nabla d)={\bf n},\qquad
  \partial_{\bf n}(\nabla d)=0.
\eqno{(B.1)}
$$
  Extend the normal vector field on $\Bbb S^2$ into $\Bbb R^3$ by setting
  ${\bf n}=\nabla d$. From the proof of Lemma 2.3 we know that there exists a
  unique $(\bfv,p,\zeta)\in (h^{m-1+\theta}(\Bbb S^2))^3\times
  h^{m-2+\theta}(\Bbb S^2)\times\Bbb R^6$, where $\bfv=(v_1,v_2,v_3)$,
  such that
$$
  \mathcal L(0)(\bfv,p,\zeta)=(0,0,\bfh,0,0).
\eqno{(B.2)}
$$
  We have
$$
  \Delta p=\nabla\cdot(\Delta\bfv-l_0(\zeta))=0
  \qquad\mbox{in}\;\;\Bbb B^3.
\eqno{(B.3)}
$$
  Since, on one hand,
\begin{eqnarray*}
  \bfT(\bfv,p)\bfn\cdot\bfn\;&=&\;2\bfn_i\bfn_j
  \partial_i v_j-p=2\bfn_j\partial_\bfn v_j-p
  \\
  &=&\;2\partial_\bfn(\bfv\cdot\bfn)-2\bfv\cdot
  \partial_\bfn(\nabla d)-p
  \\
  &=&\;2\partial_\bfn(\bfv\cdot\bfn)-p
  \qquad\mbox{on}\;\;\Bbb S^2,
\end{eqnarray*}
  and, on the other hand, by $(B.2)$,
\begin{eqnarray*}
  \bfT(\bfv,p)\bfn\cdot\bfn=\;\bfh\cdot\bfn=0
  \qquad\mbox{on}\;\;\Bbb S^2.
\end{eqnarray*}
  Hence,
$$
  2\partial_\bfn(\bfv\cdot\bfn)-p=0\qquad\mbox{on}
  \;\;\Bbb S^2.
\eqno{(B.4)}
$$
  Define $\Psi=\bfv\cdot\bfn-{1\over2}p\,d$. Then we have
  $\mbox{\rm tr}_{\Bbb S^2}(\bfv)\cdot{\bf n}=\mbox{\rm tr}_{\Bbb S^2}(\Psi)$ and
$$
\begin{array}{lll}
  \Delta \Psi\;&=&\;\Delta(\bfv\cdot\bfn-{1\over2}p\,d)
  \nonumber\\
  &=&\;\Delta \bfv\cdot\bfn+2\partial_iv_j\partial_i
  \bfn_j+\bfv\cdot\Delta\bfn-{1\over2}p\Delta d-
  \nabla d\cdot\nabla p-{1\over2}d\Delta p
  \\
  &=&\;(\nabla p+l_0(\zeta))\cdot\bfn
  +2\partial_iv_j\partial_i\bfn_j+\bfv\cdot\Delta\bfn
  -{1\over2}p\Delta d-\nabla p\cdot\bfn
  \\
  &=&\;l_0(\zeta)\cdot\bfn+2\partial_iv_j\partial_i
  \bfn_j+\bfv\cdot\Delta\bfn-{1\over2}p\Delta d
  \qquad\mbox{in}\;\;\Bbb B^3.
\end{array}
$$
  Furthermore, by $(B.1)$ and $(B.4)$ we have
$$
  \partial_\bfn\Psi=\partial_\bfn(\bfv\cdot\bfn)-{d\over2}
  \partial_\bfn p-{p\over2}\bfn\cdot(\nabla d)=\partial_\bfn
  (\bfv\cdot\bfn)-{1\over2}p=0 \qquad \mbox{on}\;\;\Bbb S^2.
$$
  Hence $\Psi$ is the solution of the problem
$$
  \left\{
\begin{array}{l}
  \Delta \Psi=\;l_0(\zeta)\cdot\bfn+2\partial_iv_j\partial_i
  \bfn_j+\bfv\cdot\Delta\bfn-{1\over2}p\Delta d
  \qquad\mbox{in}\;\;\Bbb B^3,
  \\
  \partial_\bfn\Psi=0 \qquad \mbox{on}\;\;\Bbb S^2.
\end{array}
  \right.
\eqno{(B.5)}
$$
  It follows by classical H\"older estimates for second-order partial
  differential equations of the elliptic type that
\begin{eqnarray*}
  \|\mbox{tr}_{\Bbb S^2}(\bfv)\cdot\bfn\|_{h^{m+\theta}
  (\Bbb S^2)}\;&=&\;\|\Psi\|_{h^{m+\theta}(\Bbb S^2)}\le
  C\|\Psi\|_{h^{m+\theta}(\overline{\Bbb B}^3)} \\
  &\le& C(\|\Delta \Psi\|_{h^{m-2+\theta}(\overline{\Bbb B}^3)}
  +\|\partial_{\bf n}\Psi\|_{h^{m-1+\theta}(\Bbb S^2)}+
  \|\Psi\|_{h^{m-2+\theta}(\overline{\Bbb B}^3)})
  \\
  &\le& C(\|\bfv\|_{h^{m-1+\theta}(\overline{\Bbb B}^3)}+
  \|p\|_{h^{m-2+\theta}(\overline{\Bbb B}^3)}+|\zeta|)
  \\
  &\le& C\|\bfh\|_{(h^{m-2+\theta}(\Bbb S^2))^3}.
\end{eqnarray*}
  The proof is complete. $\qquad$$\Box$
\medskip

  Since $\nabla_\omega(\partial_{\bfn}\psi_\eta)\cdot\bfn=\omega\cdot
  \nabla_\omega(\partial_{\bfn}\psi_\eta)=0$, by Lemma B.1 we immediately obtain
\medskip

  {\bf Corollary B.2}\ \ {\em Let $B_1\eta=-2{\bf R}(0)(\nabla_\omega(\partial_{\bfn}
  \psi_\eta))\big|_{\Bbb S^2}\cdot\bfn$ for $\eta\in h^{m+\theta}(\Bbb S^2)$.
  Then we have $B_1\in L(h^{m+\theta}({\Bbb S}^2),h^{m+\theta}({\Bbb S}^2))$.}
  $\qquad$ $\Box$
\medskip

   {\bf Acknowledgement}.\hskip 1em This work is supported by the National
   Natural Science Foundation of China under the grant number 10771223.

\vsss
\end{document}